\documentclass[12pt,A4]{article}
\usepackage[cp1251]{inputenc}
\usepackage{srcltx}
\usepackage{amsfonts}
\usepackage{amsmath}

\usepackage{eucal}

\begin{document}

\newcounter{abcd}
\setcounter{abcd}{1}

\fontsize{14.4}{18pt}\selectfont

\noindent UDC 517.9 \vspace{6.0mm}

\noindent \textbf{A. M. Samoilenko}\quad (Institute of Mathematics, Ukrainian
National Academy of Sciences, Kyiv)\vspace{6.0mm}

\noindent \textbf{ON INVARIANT MANIFOLDS OF LINEAR DIFFERENTIAL
EQUATIONS.~II} \vspace{8.0mm}

We continue the investigations begun in  [\,1\,]. \vspace{6.0mm}

\begin{center}\textbf{3. Equivalence of Linear Differential Equations of Different Orders}
\end{center}

As in [\,1\,], let $\Phi_{1}(t)$ and $\Phi_{2}(t)$ denote blocks of a
nonsingular matrix
\begin{equation}
\Phi(t)=\left(\begin{array}{cc}\Phi_{1}(t)\\\Phi_{2}(t)\end{array}\right)
 \nonumber
\end{equation}
that is continuously differentiable for all  $t\in \mathbb{R},$  let $\Phi_{1}^+ (t)$
and $\Phi_{2}^+ (t)$ denote blocks of the matrix
\begin{equation}
\Phi^{-1}(t)=(\Phi_{1}^+ (t), \Phi_{2}^+ (t)), \nonumber
\end{equation}
\begin{equation}
\Phi_{1}(t)\in \textbf{M}_{n \,m}(\mathbb{R}),\qquad \Phi_{2}(t)\in
\textbf{M}_{m-n \,m}(\mathbb{R}), \nonumber
\end{equation}
\begin{equation}
\Phi_{1}^+(t)\in \textbf{M}_{m \,n}(\mathbb{R}),\qquad
\Phi_{2}^+(t)\in \textbf{M}_{m \,m-n}(\mathbb{R}), \nonumber
\end{equation}
inverse to $\Phi(t),$ let $M_{1}(t)=\Phi_{1}^{+}(t)\Phi_{1}(t)$ and
$M_{2}(t)=\Phi_{2}^{+}(t)\Phi_{2}(t)$ denote projectors of ranks  $n$ and
$m-n,$ respectively, let  $M^n(t)$ and $M^{m-n}(t)$ denote the hyperplanes
\begin{equation}
 M^{n}(t)=\{y\in \mathbb{R}^{m}: y=M_{1}(t)y\},  \nonumber
\end{equation}
\begin{equation}
 M^{m-n}(t)=\{y\in \mathbb{R}^{m-n}: M_{1}(t)y=0\}  \nonumber
\end{equation}
of dimensions  $n$ and $m-n,$ respectively, and let  $L(M, Q)$ denote a matrix
operator of the form
\begin{equation}
L(M, Q)=\frac{dM}{dt} +MQ-QM. \nonumber
\end{equation}

\textbf{Theorem 2.} {\em If the subspaces  \/} $M^{n}(t) $ {\em and \/}
$M^{m-n}(t),$ {\em taken together, are invariant manifolds of the differential
equation}
\begin{equation}  \label{rim1}
\frac{dy}{dt}=Q(t)y, \tag{\Roman{abcd}}
\end{equation}
\addtocounter{abcd}{1}
{\em then the change of variables  \/}
\begin{equation}  \label{rim2}
y=\Phi_{1}^+(t)x+\Phi_{2}^+(t)z \tag{\Roman{abcd}}
\end{equation}
\addtocounter{abcd}{1} {\em reduces this equation to the system of differential
equations \/}
\begin{equation}  \label{rim3}
\frac{dx}{dt}=P(t)x, \qquad \frac{dz}{dt}=G(t)z \tag{\Roman{abcd}}
\end{equation}
\addtocounter{abcd}{1} {\em with coefficient matrices}
\begin{equation}\label{rim4}
P(t)= \left (\frac{d\Phi_1(t)}{dt}+\Phi_1(t)Q(t)
\right)\Phi_1^{+}(t), \tag{\Roman{abcd}}
\end{equation}
\addtocounter{abcd}{1}
\begin{equation}\label{rim5}
G(t)= \left (\frac{d\Phi_2(t)}{dt}+\Phi_2(t)Q(t)
\right)\Phi_2^{+}(t), \tag{\Roman{abcd}}
\end{equation}
{\em and vice versa, if the differential equation } (\textrm{I}) {\em can be reduced by the change of variables} (\textrm{II}) {\em to the system of
differential equations} (\textrm{III}), {\em then the subspaces  \/} $M^{n}(t)
$ {\em and \/}  $M^{m-n}(t),$ {\em taken together, are invariant manifolds of
Eq.}~(\textrm{I}), {\em and the coefficient matrices of system }
(\textrm{III}) {\em are defined by relations } (\textrm{IV}) {\em and }
(\textrm{V}).
\addtocounter{abcd}{1}

Indeed, let the subspaces  $M^{n}(t) $  and $M^{m-n}(t),$ taken together, be
invariant manifolds of Eq.~(I). Then, according to assertion 2 of the main
theorem in [\,1\,], Eq.~(I) is equivalent on $M^{n}(t) $ and $M^{m-n}(t)$ to
the corresponding first and second equations of system (III) with
the coefficient matrices defined by relations (IV) and (V),
respectively. Let  $Y(t),$ $X(t), $ and $Z(t),$ where $Y(0)=E,$ $X(0)=E,$ and
$Z(0)=E,$ be the fundamental matrices of solutions of Eqs.~(I) and (III) and let
$E$ be the identity matrices of the corresponding orders. According to assertion 2
of the main theorem in [\,1\,], we have
\begin{equation}  \label{1}
Y(t)\Phi^{+}_1(0)=\Phi^{+}_1(t)X(t), \qquad
Y(t)\Phi^{+}_2(0)=\Phi^{+}_2(t)Z(t)
\end{equation}
for all  $t\in \mathbb{R}.$ Thus, according to (\ref{1}),
\begin{equation}  \label{2}
Y(t)(\Phi^{+}_1(0), \Phi^{+}_2(0))=
(\Phi^{+}_1(t),\Phi^{+}_2(t))\left(\begin{array}{cc} X(t)&0
\\ 0 & Z(t)\end{array}\right)
\end{equation}
for all  $t\in \mathbb{R}.$ The equality
\begin{equation}
\label{3}Y(t)=\Phi^{+}_1(t)X(t)\Phi_1(0) +
\Phi^{+}_2(t)Z(t)\Phi_2(0)
\end{equation}
for all $t\in \mathbb{R}$ follows from  (\ref{2}). Thus, for an arbitrary $y_0
\in \mathbb{R}^m,$ we have
\begin{equation}  \label{4}
Y(t)y_0=\Phi^{+}_1(t)X(t)x_0 + \Phi^{+}_2(t)Z(t)z_0
\end{equation}
for all $t\in \mathbb{R}$ and $x_0$ and $z_0$ chosen according to the condition
\begin{equation}
x_0=\Phi_1(0)y_0, \qquad z_0=\Phi_2(0)y_0.  \nonumber
\end{equation}
Equality (\ref{4}) means that the change of variables (II) reduces the
differential equation (I) to the system of differential equations (III).

Now assume that the differential
equation (I) can be reduced to the system of differential equations (III) by the change of variables (II). Taking into account that
the subspaces $z=0$ and $x=0$ are invariant manifolds of system (III) and using (II), we obtain relations (\ref{1}), which
yield
\begin{equation}  \label{5}
X(t)=\Phi_1(t)Y(t)\Phi_1^{+}(0), \qquad
Z(t)=\Phi_2(t)Y(t)\Phi_2^{+}(0)
\end{equation}
for all $t\in \mathbb{R}.$

Substituting (\ref{5}) into relations (\ref{1}), we obtain
\begin{equation}  \label{6}
Y(t)\Phi^{+}_1(0)=M_1(t)Y(t)\Phi^{+}_1(0), \qquad
Y(t)\Phi^{+}_2(0)=M_2(t)Y(t)\Phi^{+}_2(0)
\end{equation}
for all $t\in \mathbb{R}.$

It follows from the first relation in (\ref{6}) that
\begin{equation}  \label{7}
y(t)=M_1(t)y(t)
\end{equation}
for any solution  $y(t)$ of Eq.~(I) that satisfies the condition
\begin{equation}  \label{8}
y(0)=\Phi^{+}_1(0)c,
\end{equation}
where $c$ is an arbitrary constant from  $\mathbb{R}^n.$ Since points (\ref{8})
fill the subspace  $M^n(0),$ we conclude that, according to (\ref{7}), the integral curves $(t,
y(t))$ of Eq.~(I) that pass through points of the subspace
$M^n(0)$ for  $t=0$ belong to the subspace  $M^n(t)$ for any  $t\in \mathbb{R}.$ This is
sufficient for the subspace  $M^n(t)$ to be an invariant manifold of Eq.~(I).

It follows from the second relation in (\ref{6}) that
\begin{equation}  \label{9}
y(t)=M_2(t)y(t)
\end{equation}
for any solution  $y(t)$ of Eq.~(I) that satisfies the condition
\begin{equation}  \label{10}
y(0)=\Phi^{+}_2(0)c,
\end{equation}
where $c$ is an arbitrary constant from $\mathbb{R}^{m-n}.$

By analogy, we prove that the subspace
\begin{equation}
M_2^{m-n}(t)=\{y \in \mathbb{R}^{m-n}: \; y=M_2(t)y \} \nonumber
\end{equation}
is an invariant manifold of Eq.~(I).

According to Lemma~1 in [\,1\,], the equality
\begin{equation}
M_2^{m-n}(t)= M^{m-n}(t) \nonumber
\end{equation}
holds for any  $t\in \mathbb{R}.$ This proves that the subspace $M^{m-n}(t)$
is an invariant manifold of the differential equation (I). Thus, the subspaces
$M^{n}(t)$ and $M^{m-n}(t),$ taken together, are invariant manifolds of Eq.~(I).
According to assertion 2 of the main theorem in [\,1\,], this is sufficient for
relations (IV) and (V) to be true.

Let $F(t) \in M_{p \;n}(\mathbb{R}),$ $n>p,$ $F^+(t) \in M_{p
\;n}(\mathbb{R}),$ $\textrm{rank} F(t)=p,$ and let $F(t)$ and $F^+(t)$ be continuously
differentiable functions for all $t\in \mathbb{R}.$ Also assume that $F^+(t)$
is a matrix pseudoinverse to the matrix  $F(t)$ and $K(t)=F^+(t)F(t).$ Finally, let
the subspace
\begin{equation}
K^p(t)=\{ x \in \mathbb{R}^{n}: \; x=K(t)x \} \nonumber
\end{equation}
be an invariant manifold of the differential equation
\begin{equation}  \label{rim6}
\frac{dx}{dt}=P(t)x, \tag{\Roman{abcd}}
\end{equation}
 \addtocounter{abcd}{1}
which is equivalent on  $K^p(t)$ to the differential equation
\begin{equation}  \label{rim7}
\frac{dz}{dt}=R(t)z. \tag{\Roman{abcd}}
\end{equation}
\addtocounter{abcd}{1}

The system of differential equations (III) is called a \textit{decomposition} of the
differential equation (I) if the change of variables (II) reduces this equation
to the system of differential equations (III).

The differential equation (VII) is called a \textit{restriction} of the differential
equation (VI) to the subspace  $K^p(t)$ if the subspace $K^p(t)$ is an
invariant manifold of Eq.~(VI), and this equation is equivalent to
Eq.~(VII) on $K^p(t).$

We say that the differential equations (I) and (VI) are equivalent if Eq.~(VI),
together with its restriction to  $K^{m-n}(t)$ (VII), is a decomposition of
Eq.~(I).

By definition, the fundamental matrices of solutions of equivalent
differential equations are expressed in terms of one another via the matrices that
define the invariant subspaces of these differential equations. Indeed, using the definitions presented above and taking into account that
\begin{equation}  \label{11}
G(t)=R(t)
\end{equation}
for all $t\in \mathbb{R},$ we conclude that relation (\ref{1}) and the relation
\begin{equation}  \label{12}
X(t)F^+(0)=F^+(t)Z(t)
\end{equation}
for the fundamental matrices of the solutions  $Y(t),$ $X(t),$ and $Z(t)$
of the differential equations (I), (VI), and (VII) are true.

It follows from (\ref{1}) and (\ref{12}) that
\begin{equation}
Y(t)=(\Phi_{1}^{+}(t)X(t)+\Phi_{2}^{+}(t)F(t)X(t)F^{+}(0))
\left(\begin{array}{cc}\Phi_{1}(0)\\\Phi_{2}(0)\end{array}\right)=\nonumber
\end{equation}
\begin{equation}  \label{13}
=\Phi_{1}^{+}(t)X(t)\Phi_{1}(0)+\Phi_{2}^{+}(t)F(t)X(t)F^{+}(0)\Phi_{2}(0),
\end{equation}
\begin{equation}  \label{14}
X(t)=\Phi_{1}(t)Y(t)\Phi_{1}^{+}(0)
\end{equation}
for all $t\in \mathbb{R} .$ Relations (\ref{13}) and (\ref{14}) describe the
relationship between the fundamental matrices of solutions of the equivalent
differential equations (I) and (VI).

The notion of equivalence of differential equations of orders  $m$ and $n$
defined above for
\begin{equation}
m>n>m-n \nonumber
\end{equation}
can easily be generalized to the case
\begin{equation}  \label{15}
m=2n.
\end{equation}
Indeed, since the space  $\mathbb{R}^{n}$ is an invariant manifold of the
differential equation (VI), and Eq.~(VI) is equivalent on it to the
differential equation (VII) with the same coefficient matrix, we conclude that, in case
(\ref{15}), the equivalence of the differential equations (I) and (VI) is determined
by the decomposition of Eq.~(I) into the system of equations
\begin{equation}
\frac{dx}{dt}=P(t)x,\qquad \frac{dz}{dt}=P(t)z. \nonumber
\end{equation}

The results presented above yield the following statement:

\textbf{Corollary.} \emph{The differential equations } (I) \emph{and} (VI)
\emph{are equivalent if and only if}
$$
L(M(t), Q(t))=0, \qquad L(K(t), P(t))K(t)=0,  \eqno (\textrm{VIII})
$$
$$
P(t)=\left (\frac{d \Phi_1(t)}{dt}+\Phi_1(t)Q(t) \right )\Phi_1^+(t),  \eqno (\textrm{XI})
$$
$$
\left (\frac{d \Phi_2(t)}{dt}+\Phi_2(t)Q(t) \right )\Phi_2^+(t)=
 \left (\frac{d F(t)}{dt}+F(t)P(t) \right )F^+(t)  \eqno (\textrm{X})
 $$
\emph{for all} $t \in\mathbb{R}.$

Indeed, assume that the differential equations (I) and (VI) are equivalent. Then we have the
decomposition of Eq.~(I) into the system of equations (III) the second equation
of which is the restriction of the differential equation (VI) to $K^{m-n}(t).$ It follows from the definition of decomposition and Theorem~2 that
the subspaces $M^{n}(t)$ and $M^{m-n}(t)$ are invariant manifolds of the
differential equation (I). It follows from the definition of the restriction of
the differential equation (VI) to the subspace  $K^{m-n}(t)$ that $K^{m-n}(t)$
is an invariant manifold of this equation. According to assertions 1 and 2 of
the main theorem in [\,1\,], this is sufficient for relations (VIII) and (IX) to be
true. Moreover, this is sufficient for the coefficient matrices of the
differential equations (I), (III), and (VII) to satisfy the relations
\begin{equation}  \label{16}
 G(t)= \left (\frac{d \Phi_2(t)}{dt}+ \Phi_2(t)Q(t) \right )
 \Phi_2^+(t),
\end{equation}
\begin{equation}  \label{17}
R(t)=\left (\frac{d F(t)}{dt}+F(t)P(t) \right )F^+(t),
\end{equation}
and
\begin{equation}  \label{18}
G(t)=R(t)
\end{equation}
for all  $t \in\mathbb{R}.$

The last relation proves equality (X).

Let relations (VIII) -- (X) be true. Then, according to assertions 1 and 2 of
the main theorem in [\,1\,], the subspaces  $M^{n}(t)$ and $M^{m-n}(t)$ are
invariant manifolds of the differential equation (I), and the subspace
$K^{m-n}(t)$ is an invariant manifold of the differential equation (VI);
furthermore, the coefficient matrices of the corresponding differential
equations $G(t)$ and $R(t)$ are defined by relations (\ref{16}) and (\ref{17}),
and, hence, according to condition (X), they satisfy equality (\ref{18}). According to
Theorem~2, this implies that the system of differential equations (III) the
second equation of which is the restriction of the differential equation (VI) to
the subspace  $K^{m-n}(t)$ is a decomposition of the differential equation
(I). This proves that relations (VIII)--(X) yield the equivalence of the
differential equations (I) and (VI).

Note that, for $m=2n,$ conditions (VIII)--(X) are simplified because, in this case, $F(t)$ and $K(t)$ are the identity matrices. In this case, these conditions take the form
\begin{equation}
 L(M(t), Q(t))=0, \qquad
  P(t)=\left (\frac{d \Phi_1(t)}{dt}+\Phi_1(t)Q(t)
\right )\Phi_1^+(t) = \nonumber
\end{equation}
\begin{equation}= \left (\frac{d
\Phi_2(t)}{dt}+\Phi_2(t)Q(t) \right )\Phi_2^+(t)  \nonumber
\end{equation}
for any $t \in\mathbb{R}.$

Also note that the equivalence of the differential equations (I) and (VI) means
that the relations
\begin{equation}  \label{19}
Y(t) \Phi_{1}^{+}(0)=\Phi_{1}^{+}(t) X(t), \qquad Y(t)
\Phi_{2}^{+}(0)= \Phi_{2}^{+}(t) F(t) X(t) F^{+}(0)
\end{equation}
for the fundamental matrices of solutions of Eqs.~(I) and (VI) $Y(t)$ and
$X(t),$ as well as the other relations that can be obtained from (\ref{19}) by the corresponding transformations, are true.

\begin{center}\textbf{4. Addition to the Floquet--Lyapunov Theory}\end{center}

Consider the linear differential equation
 \setcounter{abcd}{1}
\begin{equation}  \label{rim1}
\frac{dx}{dt}=P(t)x, \tag{\Roman{abcd}}
\end{equation}
where  $ x \in\mathbb{R}^n, \, P(t)\in \textbf{M}_n(\mathbb{R}),$ and $P(t)$ is a continuous periodic matrix with period  $T.$
 \addtocounter{abcd}{1}

According to the well-known Floquet theorem [\:2\,], the fundamental matrix of
solutions of Eq.~(I) $X(t), \, X(0)=E, $ can be represented in the form
$$
X(t)=\Phi (t) e^{Ht},  \eqno(\textrm{II})
$$
where $\Phi(t)$ is a matrix periodic in $t$ with period $T,$ and
$H$ is the constant matrix defined by the monodromy matrix  $X(T)$ of Eq.~(I) according to
the formula
$$
H=\frac{1}{T}\ln X(T).  \eqno(\textrm{III})
$$
The logarithm is a multi-valued function whose real value does not always exist. Thus,
relation (I) with matrix (III) such that
$$
H \in \textbf{M}_n(\mathbb{R})  \eqno(\textrm{IV})
$$
is not always true. According to the theory of matrices [\,3\,], condition (IV)
is satisfied if and only if every elementary divisor corresponding to the negative
eigenvalues of the matrix  $X(T)$ is repeated an even number of times. Thus,
only in this case does equality (II) hold with matrices  $\Phi(t)$ and $H$ from
the space of real matrices $\textbf{M}_n(\mathbb{R}).$

In the case where condition (IV) cannot be satisfied, the Floquet representation
(II) exists only with matrices  $\Phi(t)$ and $H$  from the space
$\textbf{M}_n(\mathbb{C}),$ where $\mathbb{C}$ is the plane of complex numbers,
or this representation transforms into equality (II) with real matrices
$\Phi(t)$ and $H,$ the first of which is periodic with period  $2T$ and the
second is defined by the relation
$$
H=\frac{1}{2T}\ln X(2T).  \eqno(\textrm{V})
$$
The Floquet representation (II) with matrix (V) is a consequence of the presence
of negative eigenvalues of the monodromy matrix of Eq.~(I).

We consider in detail the differential equation (I) whose
monodromy matrix possesses this property and prove several previously unknown statements for this equation.

\textbf{Theorem 3.} {\em Suppose that the coefficient matrix of the differential
equation} (I)  $ P(t)$ \textit{belongs to} $\textbf{M}_n(\mathbb{R})$ {\em for any $t
\in\mathbb{R}$ and is continuous on  $\mathbb{R}$ and periodic in  $t $ with
period  $T.$ }

{\em Then the following assertions are true:\/}

1. {\em The algebraic number $p$ of negative eigenvalues of the monodromy matrix
$X(T)$ of Eq.~}(I) {\em is even. \/}

2. {\em Equality} (II) {\em holds for the matrix \/}
$$
H=\frac{1}{T}\ln (X(T)I),  \eqno(\textrm{VI})
$$
{\em  where $I$ is the real matrix defined by the conditions \/}
$$
I^2=E, \qquad \ln (X(T)I) \in \textbf{M}_n(\mathbb{R}),
$$
{\em and for the periodic matrix  $\Phi(t)$ such that \/}
$$
\Phi(t+T)I_1=\Phi(t)I_1, \qquad   \Phi(t+T)I_2=-\Phi(t)I_2  \eqno(\textrm{VII})
$$
{\em for all  $t \in\mathbb{R},$  where \/}
$$
I_1=\frac{E+I}{2}, \qquad I_2=\frac{E-I}{2}.
$$

3. {\em There exists a nonsingular matrix $(U(t), \, V(t))$ continuously
differentiable and real for all  $t \in\mathbb{R},$ periodic with period $T,$ and
such that the change of variables
$$
x= U(t)z_1+ V(t)z_2
$$
reduces the differential equation} (I) {\em to the system of differential
equations \/}
$$
\frac{d z_1}{dt}=H_1 z_1, \qquad \frac{d z_2}{dt}=G(t) z_2,  \eqno(\textrm{VIII})
$$
{\em where $H_1$ is a constant matrix, $G(t)$ is a periodic matrix with period
$T,$ and the set of eigenvalues of the monodromy matrix  $Z_2(T)$ of the second
equation of the system is either the set of all negative eigenvalues of the
matrix $X(T)$ or its subset. \/}

To prove the theorem, we use the representation of the matrix  $X(T)$
in terms of its Jordan form $J(\lambda),$ namely
\begin{equation}
 X(T)=S J(\lambda)S^{-1},\nonumber
 \end{equation}
and obtain the equality
 \setcounter{equation}{0}
\begin{equation}  \label{1}
\det X(T)=\prod ^n_{\nu =1} \lambda _ \nu,
\end{equation}
which associates the determinant of the matrix $X(T)$ with its eigenvalues
$\lambda _ \nu\/, \, \nu=\overline{1,\, n}.$

We now use the Liouville--Ostrogradskii--Jacobi formula and represent the
determinant of the matrix  $X(T)$ in terms of the trace of the coefficient matrix
of Eq.~(I):
\begin{equation}  \label{2}
 \det X(T)= \exp \left \{ \textrm{tr} \,P(t) dt \right \} .
\end{equation}
Equating the right-hand sides of relations (\ref{1}) and (\ref{2}), we obtain
an equality that proves that
\begin{equation}  \label{3}
\prod ^n_{\nu =1} \lambda _ \nu > 0.
\end{equation}
Since each pair of complex conjugate eigenvalues of the matrix  $X(T)$ in the
product of all its eigenvalues gives a positive number, it follows from relation (\ref{3})
that the product of all negative eigenvalues of the matrix  $X(T)$
also gives a positive number. Thus, the algebraic number of negative eigenvalues of
the matrix  $X(T),$ i.e., the sum of multiplicities of the roots of characteristic equations for all different negative eigenvalues of
the matrix $X(T),$ is an even
number.

Prior to the proof of assertion 2 of Theorem~3, note that, in the case where
the logarithm of the matrix  $X(T)$ is real, by setting $I=E$ one can reduce
equalities (II) and (VI) to the Floquet relations (II) and (III) with a matrix
$\Phi(T)$ that possesses properties that follow from these relations and are indicated in
assertion 2 of Theorem~3.

It remains to consider the case where the matrix  $X(T)$ has negative
eigenvalues and does not have a real logarithm. In this case, the real
canonical form of the matrix  $X(T)$ can be represented in the form of decomposition
into two blocks  $A$ and $B,$ where  $A$ either is empty or has a real
logarithm, and  $B$ has only negative eigenvalues and does not have a real
logarithm.

Let $B\in \textbf{M}_d(\mathbb{R}),$ where
\begin{equation}  \label{4}
n>d.
\end{equation}
Then the following equality is true:
\begin{equation}  \label{5}
X(T)=S \left(\begin{array}{cc} A&0
\\ 0 & B\end{array}\right) S^{-1}\, ,
\end{equation}
where $S,$ $A,$ and $B$ are real matrices with properties indicated above for
$A$ and~$B.$

We set
\begin{equation}  \label{6}
Y(t)=S^{-1}X(t)S,\qquad B_{1}=-B.
\end{equation}
According to properties of the fundamental matrix of solutions of Eq.~(I), we have
\begin{equation}  \label{7}
X(t+T)=X(t)X(T)\,.
\end{equation}
Therefore, it follows from (\ref{5}), (\ref{6}), and (\ref{7}) that
\begin{equation}  \label{8}
Y(t+kT)=S^{-1}X(t)X^{k}(T)S=S^{-1}X(t)SS^{-1}X^{k}(T)S=\nonumber
\end{equation}
\begin{equation} =Y(t) \left(\begin{array}{cc} A^k&0
\\ 0 & B^k\end{array}\right)=Y(t)\left(\begin{array}{cc} A^k&0
\\ 0 & (-1)^k B_1^k\end{array}\right)
\end{equation}
for any integer $k.$

We represent  $Y(t)$ in the block form
\begin{equation}\label{9}
Y(t)=(Y_{1}(t),Y_{2}(t))
\end{equation}
consistent with decomposition (5) of the matrix  $X(T)$ into the blocks $A$ and $B.$
Using relations (\ref{8}), we get
\begin{equation}\label{10}
Y_{1}(t+kT)=Y_{1}(t)A^{k},\qquad
Y_{2}(t+kT)=(-1)^{k}Y_{2}(t)B_{1}^{k}
\end{equation}
for any integer $k$.

Since the eigenvalues of the matrix  $B_{1}$ are positive by virtue of the definition
(\ref{6}) of this matrix, both matrices  $A$ and $B_{1}$ have real logarithms
$\ln A$ and $\ln B_{1}.$

In view of the arguments presented above, relation (\ref{10}) yields
\begin{equation}
Y_{1}(t)=Y_{1}\left(t-\left[\frac{t}{T}\right]T+\left[\frac{t}{T}\right]T\right)=Y_{1}\left
(t- \left[\frac{t}{T}\right]T \right)A^{\left[
\frac{t}{T}\right]}=\nonumber
\end{equation}
\begin{equation}\label{11}
=Y_{1}\left(t- \left[\frac{t}{T}\right]T \right)\exp \left \{
\left(\left[ \frac{t}{T}\right]T-t\right) \frac{\ln A}{T} \right \}
 \exp \left \{ \frac{t}{T} \ln A \right \}  ,
\end{equation}
\begin{equation}
Y_{2}(t)=Y_{2}\left(t-\left[\frac{t}{T}\right]T+\left[\frac{t}{T}\right]T\right)=Y_{2}\left(t-
\left[\frac{t}{T}\right]T\right)(-1)^{\left[\frac{t}{T}\right]}B_{1}^{\left[\frac{t}{T}\right]}=\nonumber
\end{equation}
\begin{equation}\label{12}
=(-1)^{\left[\frac{t}{T}\right]}Y_{2}\left(t-
\left[\frac{t}{T}\right]T\right) \exp \left \{
\left(\left[\frac{t}{T}\right]T-t\right)\frac{\ln B_{1}}{T} \right
\} \exp \left \{\frac{t}{T}\ln B_{1} \right \}
\end{equation}
for all  $t\in\mathbb{R};$ here, $[\,t\,]$ denotes the integer part of the
number  $t.$

Let $\Phi_{1}(t)$ and $\Phi_{2}(t)$ denote the coefficients of $\, \exp \left
\{\displaystyle \frac{t}{T}\ln A \right \}$ and $\, \exp \left \{\displaystyle
\frac{t}{T}\ln B_{1} \right \}$ in relations (\ref{11}) and (\ref{12}),
respectively. Then, using (\ref{9}), (\ref{11}), and (\ref{12}), we obtain
\begin{equation}\label{13}
Y(t)=(\Phi_{1}(t),\Phi_{2}(t)) \left(\begin{array}{cc} \exp \left
\{\displaystyle \frac{t}{T}\ln A \right \} &0
\\ 0 & \exp \left \{\displaystyle \frac{t}{T}\ln B_{1} \right \} \end{array}\right)
\end{equation}
for all $t\in\mathbb{R}.$ This equality implies that the matrices
$\Phi_{1}(t)$ and $\Phi_{2}(t)$ are continuously differentiable on
$\mathbb{R}.$ Furthermore, it follows from the introduced notation that the matrix
$\Phi_{1}(t)$ is periodic with period  $T,$ and the matrix $\Phi_{2}(t),$ which
is the product of the function  $(-1)^{[ \frac{t}{T}]}$ and a periodic matrix with
period  $T,$ satisfies the condition
\begin{equation}\label{14}
\Phi_{2}(t+T)=-\Phi_{2}(t)
\end{equation}
for all $t\in\mathbb{R}.$

Let $I_{0}$ denote the matrix
$$
\left( \begin{array}{cc} E_1 & 0 \\ 0 & -E_2 \end{array} \right),
$$
where $E_{1}$ and $E_{2}$ are the identity matrices from
$\textbf{M}_{n-d}(\mathbb{R})$ and $\textbf{M}_{d}(\mathbb{R}),$ respectively.
Then
$$
Y(T)I_{0}=\left( \begin{array}{cc} A & 0 \\ 0 & B_1 \end{array} \right)
$$
and relation (\ref{13}) takes the form
\begin{equation}\label{15}
Y(t)=(\Phi_{1}(t),\Phi_{2}(t)) \exp \left \{ \frac{t}{T}\ln (Y(T)I_{0}) \right \}.
\end{equation}
Using (\ref{15}) and the first equality in (\ref{6}), we obtain
\begin{equation}\label{16}
X(t)=S(\Phi_{1}(t),\Phi_{2}(t))S^{-1}\exp \left \{ \frac{t}{T}S(\ln(Y(T)I_{0}))S^{-1} \right \} .
\end{equation}
Since
\begin{equation}\label{17}
S(\ln(Y(T)I_{0}))S^{-1}=\ln(SY(T)S^{-1}SI_{0}S^{-1})=\ln(X(T)I),
\end{equation}
where
\begin{equation}  \label{18}
I=SI_{0}S^{-1},
\end{equation}
relation (\ref{16}) takes the form of the required representation (II) under the condition
that
\begin{equation}\label{19}
H=\frac{1}{T}\ln(X(T)I),
\end{equation}
\begin{equation}\label{20}
\Phi(t)=S(\Phi_{1}(t), \Phi_{2}(t))S^{-1}.
\end{equation}
Taking (\ref{18}) into account, we get
\begin{equation}
I^{2}=E,\qquad I_{1}=S \left( \begin{array}{cc} E_1 & 0 \\ 0 & 0
\end{array} \right)S^{-1}, \qquad I_{2}=S \left( \begin{array}{cc} 0 & 0 \\ 0 &
E_2 \end{array} \right)S^{-1}, \nonumber
\end{equation}
where $E_{1}$ and $E_{2}$ are the identity matrices of the corresponding
orders.

Using the expressions for  $I_{1}$ and $I_{2},$ we obtain
\begin{equation}\label{21}
\Phi(t)I_{1}=S(\Phi_{1}(t),0)S^{-1},\qquad\Phi(t)I_{2}=S(0,\Phi_{2}(t))S^{-1}
\end{equation}
for all $t\in\mathbb{R}.$ In view of properties of the matrices
$\Phi_{1}(t)$ and $\Phi_{2}(t),$ relation
(\ref{21}) yields
\begin{equation}
\Phi(t+T)I_{1}=\Phi(t)I_{1},\qquad
\Phi(t+T)I_{2}=-\Phi(t)I_{2}\nonumber
\end{equation}
for all  $t\in\mathbb{R},$ which completes the proof of assertion 2 of
Theorem~3 in the case considered.

Let $d=n.$ In this case, we obtain the equality
\begin{equation}
X(T)=SBS^{-1} \nonumber
\end{equation}
instead of (\ref{5}), the equality
\begin{equation}
Y(t+kT)=(-1)^{k}Y(t)B_1^k \nonumber
\end{equation}
instead of (\ref{8}), and the equality
\begin{equation}
Y(t)=\Phi_{2}(t) \exp \left \{ \frac{t}{T}\ln B_{1} \right \}
\nonumber
\end{equation}
and condition
\begin{equation}
\Phi_{2}(t+T)=-\Phi_{2}(t)\nonumber
\end{equation}
for all  $t\in\mathbb{R}$ instead of (\ref{13}).

We set
\begin{equation}
I_{0}=-E. \nonumber
\end{equation}
Using the last two formulas, we obtain equality (II) of the form
\begin{equation}
Y(t)=\Phi_{2}(t) \exp \left \{ \frac{t}{T}\ln (-Y(t)) \right \}, \nonumber
\end{equation}
where
\begin{equation}
H=\frac{1}{T}\ln (-X(T)), \qquad H \in \textbf{M}_{n}\mathbb{R},
\nonumber
\end{equation}
\begin{equation}
\Phi(t)=S\Phi_{2}(t)S^{-1}, \qquad  \Phi(t+T)=-\Phi(t),\qquad \Phi(t) \in
\textbf{M}_{n}\mathbb{R}, \nonumber
\end{equation}
for all  $t\in\mathbb{R},$ which completes the proof of assertion 2 of
Theorem~3.

We now pass to the proof of assertion 3 of Theorem~3. In this assertion,
we separate two limiting cases, namely, the case where the matrix  $X(t)$ has a
real logarithm and the second case where all eigenvalues of the matrix  $X(t)$ are
negative and their elementary divisors are different.

In the first case, assertion 3 of Theorem~3 follows from the Floquet relations
(II) and (III), according to which the change of variables
\begin{equation}
x=\Phi(t)z\nonumber
\end{equation}
reduces the differential equation (I) to the differential equation
\begin{equation}
\frac{dz}{dt}=Hz\nonumber
\end{equation}
and guarantees the properties of the matrices  $H$ and $\Phi(t)$ indicated in
Theorem~3.

In the second case, assertion 3 of Theorem 3 is trivial: the change of
variables
\begin{equation}
x=z\nonumber
\end{equation}
reduces the differential equation (I) to a differential equation with the
same coefficient matrix:
\begin{equation}
\frac{dz}{dt}=P(t)z.\nonumber
\end{equation}

Associating these limiting cases with the representation of the matrix  $X(T)$ via its real canonical form (\ref{5}), we establish that the first case corresponds to
\begin{equation}
X(T)=SAS^{-1} \nonumber
\end{equation}
and the second case corresponds to
\begin{equation}
X(T)=SBS^{-1}.\nonumber
\end{equation}
Thus, the only nonlimiting case in assertion 3 of Theorem~3 is the case where
\begin{equation}
A\in \textbf{M}_{n-d}(\mathbb{R}),\qquad B\in \textbf{M}_{d}(\mathbb{R}),
\qquad n>d>1. \nonumber
\end{equation}
Assume that these conditions are satisfied. Then it follows from the proof of
assertion 2 of Theorem 3 that the matrix  $Y(t)$ associated with the matrix
$X(t)$ by relation (\ref{6}) has the form (\ref{13}). Denoting
\begin{equation}
U(t)=\Phi_{1}(t),\qquad V(t)=\Phi_{2}(t),\qquad H_{1}=\frac{\ln
A}{T}\,,\qquad H_{2}=\frac{\ln B_{1}}{T}, \nonumber
\end{equation}
we represent (\ref{13})  in the form
\begin{equation}\label{22}
Y(t)=(U(t),V(t))\left( \begin{array}{cc} e^{H_1 t} & 0 \\ 0 & e^{H_2 t} \end{array} \right).
\end{equation}
It follows from (\ref{22}) that
\begin{equation}\label{23}
Y(t)\left( \begin{array}{cc} E_1  \\ 0  \end{array} \right)=(U(t),V(t))\left( \begin{array}{cc} e^{H_1 t}  \\ 0  \end{array} \right)=U(t)e^{H_{1}t}\,,
\end{equation}
where $E_{1}$ is the identity matrix of order $n-d.$

Differentiating equality (\ref{23}) with regard for the first relation in
(\ref{6}), we get
\begin{equation}
S^{-1}P(t)SY(t)\left( \begin{array}{cc} E_1 \\ 0 \end{array} \right)=S^{-1}P(t)SU(t)e^{H_{1}t}=\frac{dU(t)}{dt}e^{H_{1}t}+U(t)H_{1}e^{H_{1}t}.
\nonumber
\end{equation}
Thus,
\begin{equation}\label{24}
\frac{dU(t)}{dt}+U(t)H_{1}=S^{-1}P(t)SU(t)
\end{equation}
for all $t\in\mathbb{R}.$

The matrix $Y(t)$ is the fundamental matrix of solutions of the differential
equation
\begin{equation}\label{25}
\frac{dy}{dt}=S^{-1}P(t)Sy.
\end{equation}

Let $W(t)\in \textbf{M}_{n \:d}(\mathbb{R})$ for all $t\in\mathbb{R}$ and let this matrix be
continuously differentiable on  $\mathbb{R},$ periodic with period $T,$ and
such that
\begin{equation}\label{26}
\det(U(t),W(t))\ne 0
\end{equation}
for all  $t\in\mathbb{R}.$

The existence of this matrix follows from the theorem on a quasiperiodic basis
in $\mathbb{R}^{n}$ presented in  [\,4\,].

In the differential equation (\ref{25}), we perform the change of variables according to the formula
\begin{equation}\label{27}
y=U(t)y_{1}+W(t)y_{2}\,.
\end{equation}
Using equality (\ref{24}), we obtain the differential equation
\begin{equation}\label{28}
U(t)\left(\frac{dy_{1}}{dt}-H_{1}y_{1}\right)+W(t)\frac{dy_{2}}{dt}=\left(S^{-1}P(t)SW(t)
-\frac{dW(t)}{dt}\right)y_{2}\,.
\end{equation}
Solving this equation with the use of the matrix
\begin{equation}\label{29}
\left( \begin{array}{cc} L_1(t) \\ L_2(t) \end{array} \right)
\end{equation}
that is inverse to the matrix  $(U(t),W(t)),$ we obtain the following system
of differential equations for  $\displaystyle \frac{dy_{1}}{dt}$ and
$\displaystyle \frac{dy_{2}}{dt}$:
\begin{equation}\label{30}
\frac{dy_{1}}{dt}=H_{1}y_{1}+L_{1}(t)\left(S^{-1}P(t)SW(t)-\frac{dW(t)}{dt}\right)y_{2}\,,
\end{equation}
\begin{equation}\label{31}
\frac{dy_{2}}{dt}=L_{2}(t)\left(S^{-1}P(t)SW(t)-\frac{dW(t)}{dt}\right)y_{2}.
\end{equation}

Since the coefficient matrix of system (\ref{30}), (\ref{31}) has a
block-triangular form, the fundamental matrix of solutions of this system is the
matrix
\begin{equation}\label{32}
\left( \begin{array}{cc} e^{H_1 t} & Y_1(t) \\ 0 & Y_2(t) \end{array} \right)
\end{equation}
the second column of which is formed by solutions of the system of
differential equations (\ref{30}), (\ref{31}) with given initial values
$y_{1}=Y_{1}(0)$  and  $y_{2}=Y_{2}(0)$ such that
\begin{equation}\label{33}
\det\,Y_{2}(0)\ne 0.
\end{equation}
In view of (\ref{27}), the matrix
\begin{equation}\label{34}
(U(t),W(t))\left( \begin{array}{cc} e^{H_1 t} & Y_1(t) \\ 0 & Y_2(t) \end{array} \right)
\end{equation}
is a fundamental matrix of
solutions of Eq.~(\ref{25}). Moreover, relation (\ref{22}) also determines a fundamental matrix of
solutions of Eq.~(\ref{25}). According to the theory of linear differential
equations, there exists the following relation between these two fundamental
matrices of solutions:
\begin{equation}\label{35}
(U(t),V(t))\left( \begin{array}{cc} e^{H_1 t} & 0 \\ 0 & e^{H_2 t} \end{array} \right)C=(U(t),W(t))\left( \begin{array}{cc} e^{H_1 t} & Y_1(t) \\ 0 & Y_2(t) \end{array} \right)
\end{equation}
for all $t\in\mathbb{R},$ where $C$ is a nonsingular constant matrix.
Substituting $t=0$ into (\ref{35}), we obtain the following algebraic equation
for the determination of the matrix $C$:
\begin{equation}\label{36}
(U(0),V(0))C=(U(0),W(0))\left( \begin{array}{cc} E_1 & Y_1(0) \\ 0 &
Y_2(0) \end{array} \right).
\end{equation}
Multiplying (\ref{36}) by the matrix  $\left( \begin{array}{cc} L_1(0) \\
L_2(0) \end{array} \right),\,$  we obtain
\begin{equation}\label{37}
\left( \begin{array}{cc} E_1 & L_1(0)V(0) \\ 0 & L_2(0)V(0) \end{array} \right)
C = \left( \begin{array}{cc} E_1 & Y_1(0) \\ 0 & Y_2(0) \end{array} \right).
\end{equation}
This equality implies that, first,
\begin{equation}\label{38}
\det(L_{2}(0),V(0))\ne 0
\end{equation}
and, second, for
\begin{equation}\label{39}
Y_{1}(0)=L_{1}(0)V(0),\qquad Y_{2}(0)=L_{2}(0)V(0),
\end{equation}
we have
\begin{equation}\label{40}
C=E.
\end{equation}
Thus, determining the solutions  $\left(
\begin{array}{cc}  Y_1(t) \\  Y_2(t) \end{array} \right)$ of the system of differential equations
(\ref{30}), (\ref{31}) with initial values (\ref{39}), we obtain the following
equality from (\ref{35}) and (\ref{40}):
\begin{equation}\label{41}
(U(t),V(t))\left( \begin{array}{cc} e^{H_1 t} & 0 \\ 0 & e^{H_2 t} \end{array} \right)=(U(t),W(t))\left( \begin{array}{cc} e^{H_1 t} & Y_1(t) \\ 0 & Y_2(t) \end{array} \right)
\end{equation}
for all $t\in\mathbb{R}$.

Multiplying (\ref{41}) by matrix (\ref{29}), we get
\begin{equation}
\left( \begin{array}{cc} E_1 & L_1(t)V(t) \\ 0 & L_2(t)V(t)\end{array} \right) \left( \begin{array}{cc} e^{H_1 t} & 0 \\ 0 & e^{H_2 t} \end{array} \right)=\left( \begin{array}{cc} e^{H_1 t} & Y_1(t) \\ 0 & Y_2(t) \end{array} \right). \nonumber
\end{equation}
Thus,
\begin{equation}\label{42}
Y_{1}(t)=L_{1}(t)V(t)e^{H_{2}t},
\end{equation}
\begin{equation}\label{43}
Y_{2}(t)=L_{2}(t)V(t)e^{H_{2}t}
\end{equation}
for all  $t\in\mathbb{R}$. Since the matrix $Y_{2}(t)$ is nonsingular, we
can determine the value of $e^{H_{2}t}$ from (\ref{43}). Substituting this value
into (\ref{42}), we establish that
\begin{equation}\label{44}
Y_{1}(t)=L_{1}(t)V(t)(L_{2}(t)V(t))^{-1}Y_{2}(t)
\end{equation}
for all  $t\in\mathbb{R}.$

We rewrite the system of differential equations (\ref{30}), (\ref{31}) in the
form of the system
\begin{equation}\label{45}
\frac{dy_{1}}{dt}=H_{1}y_{1}+R_{1}(t)y_{2},
\end{equation}
\begin{equation}\label{46}
\frac{dy_{2}}{dt}=G(t)y_{2},
\end{equation}
where
\begin{equation}
R_{1}(t)=L_{1}(t)\left(S^{-1}P(t)SW(t)-\frac{dW(t)}{dt}\right),\nonumber
\end{equation}
\begin{equation}
G(t)=L_{2}(t)\left(S^{-1}P(t)SW(t)-\frac{dW(t)}{dt}\right).\nonumber
\end{equation}
Using the matrix
\begin{equation}\label{47}
F(t)=L_{1}(t)V(t) (L_{2}(t)V(t))^{-1},
\end{equation}
we rewrite equality (\ref{44}) in the form
\begin{equation}\label{48}
Y_{1}(t)=F(t)Y_{2}(t).
\end{equation}
Differentiating equality (\ref{48}) and taking into account that the matrix
$\left(\begin{array}{cc} Y_1(t) \\  Y_2(t) \end{array} \right)$ is a block of
the fundamental matrix (\ref{32}) of solutions of the system of differential
equations (\ref{30}), (\ref{31}) [and, hence, of system (\ref{45}), (\ref{46})], we
get
\begin{equation}\label{49}
\frac{dF(t)}{dt}+F(t)G(t)=H_{1}F(t)+R_{1}(t)
\end{equation}
for all $t\in\mathbb{R}.$ Finally, performing the change of variables
\begin{equation}\label{50}
y_{1}=z_{1}+F(t)z_{2},\qquad y_{2}=z_{2},
\end{equation}
we obtain the system
\begin{equation}
\frac{dz_{1}(t)}{dt}+\frac{dF(t)}{dt}z_{2}+F(t)G(t)z_{2}=H_{1}z_{1}+H_{1}F(t)z_{2}+
R_{1}(t)z_{2},\nonumber
\end{equation}
\begin{equation}
\frac{dz_{2}(t)}{dt}=G(t)z_{2} \nonumber
\end{equation}
instead of the system differential equations (\ref{45}), (\ref{46}). By virtue
of (\ref{49}), this system takes the form
\begin{equation}\label{51}
\frac{dz_{1}(t)}{dt}=H_{1}z_{1},\qquad
\frac{dz_{2}(t)}{dt}=G(t)z_{2}.
\end{equation}
Since the second equation of system (\ref{51}) coincides (to within notation) with
Eq.~(\ref{46}), the matrix $Y_{2}(t)$ is a fundamental matrix of solutions of
the second equation of system (\ref{51}). Then, according to relation (\ref{43}), the
matrix
\begin{equation}\label{52}
L_{2}(t)V(t)e^{H_{2}t}(L_{2}(0)V(0))^{-1}
\end{equation}
is a fundamental matrix of solutions of the second equation of system
(\ref{51}) and is equal to the identity matrix for $t=0.$ Thus, the matrix
\begin{equation}\label{53}
(L_{2}(T)V(T))e^{H_{2}T}(L_{2}(0)V(0))^{-1}
\end{equation}
is the monodromy matrix of the second equation of system (\ref{51}).

By definition, the matrix $L_{2}(t)$ is periodic with period $T,$
the matrix $V(t)$ satisfies the condition
\begin{equation}\label{54}
V(t+T)=-V(t),
\end{equation}
and the matrix  $H_{2}$ has the form
\begin{equation}\label{54}
H_{2}=\frac{1}{T}\ln (-B).\nonumber
\end{equation}
Taking into account the properties of the matrices $L_{2}(t),$ $V(t),$ and
$H_{2}$ presented above, we conclude that matrix (\ref{53}) has the form
\begin{equation}
(-L_{2}(0)V(0))(-B)(L_{2}(0)V(0))^{-1}=L_{2}(0)V(0)B(L_{2}(0)V(0))^{-1}.\nonumber
\end{equation}
Thus, it follows from the results presented above and the definition of the matrix
$B$ that the set of eigenvalues of matrix (\ref{53}) is either
the set of all negative eigenvalues of the matrix $X(T)$ or its subset.

Consider the matrix  $F(t).$ The definition of this matrix
[see (\ref{47})] and the fact that the matrices  $L_{1}(t)$ and $L_{2}(t)$ are
periodic with period  $T$ and the matrix $V(t)$ satisfies condition (\ref{54}) imply that
\begin{equation}
F(t+T)=(-L_{1}(t)V(t))(-L_{2}(t)V(t))^{-1}=F(t)\nonumber
\end{equation}
for all  $t\in\mathbb{R}.$

Thus, the matrix  $F(t)$ is periodic with period  $T.$

To complete the proof of assertion 3 of Theorem~3, it remains to take into
account that the change of variables
\begin{equation}\label{55}
x=Sy
\end{equation}
transforms the differential equation (I) into the differential equation
(\ref{25}). Therefore, the superposition of changes (\ref{55}), (\ref{27}), and
(\ref{25}) transforms the differential equation (I) into the system of
differential equations (\ref{51}), and both the change of
variables and the differential equations of system (\ref{51}) themselves possess the properties
indicated in Theorem~3.

We now make several remarks on assertions 2 and 3 of Theorem~3.

The first remark deals with relation (VI), which defines the matrix $H.$ It
follows from the proof of Theorem~3 that  $H$ is not always uniquely defined.
This nonuniqueness is caused by the condition of decomposition of the canonical
form of the matrix $X(T)$ into blocks $A$ and $B$ according to which the
matrix $B$ can be either a block of the Jordan form of the matrix $X(T)$ formed
by all its Jordan cells corresponding to its negative eigenvalues or a block
of this form obtained from the block indicated above by elimination of an arbitrary
number of pairs of identical Jordan cells.

The second remark deals with the minimum possible value of the order of the second differential equation of system (VIII). It follows from the proof of Theorem~3 that this order is also related to the condition of decomposition of the real canonical form of the matrix $X(T)$ into blocks $A$ and $B$ and is equal to the minimum possible order of the matrix $B$ of this decomposition. It follows from the first remark that the minimum possible value of the order of the second equation of system (VIII) is equal to the order of the matrix obtained from the Jordan form of the matrix $X(T)$ by elimination of all Jordan cells corresponding to negative eigenvalues of the matrix $X(T)$ and the maximum possible even number of identical Jordan cells of this matrix that correspond to its negative eigenvalues.

Also note that, according to the proof of Theorem~3, the matrix $B$ is the Jordan form of the monodromy matrix $Z_2(T)$ of the second equation of system (VIII), and, hence, the fundamental matrix of solutions $Z_2(t)$ of this equation possesses all the corresponding properties.

Finally, note that, by virtue of Theorem~2 and assertion 3 of Theorem~3, the
differential equation (I) has the invariant manifolds
\begin{equation}
K^{n-d}(t)=\{ x\in \mathbb{R}^n:\, U(t)L_1(t)x=x \}, \nonumber
\end{equation}
\begin{equation}
K^{d}(t)=\{ x\in \mathbb{R}^n:\, V(t)L_2(t)x=x \} \nonumber
\end{equation}
periodic with period $T$,
\begin{equation}
K^{\nu}(t+T)=K^{\nu}(t), \qquad \nu \in \{(n-d) \vee d\}, \nonumber
\end{equation}
for all $t\in \mathbb{R} .$ Moreover, Eq.~(I) is equivalent on  $K^{n-d}(t)$
to the first differential equation of system (VIII) and on $K^{d}(t)$ to the second differential equation of this system.

\textbf{Corollary.} \emph{The fundamental matrix of solutions of the
differential equation} (I) $X(t)$ \emph {admits the representation}
$$
X(t)=\Phi (t) e^{Ht}\Phi ^+(0),  \eqno (\textrm{IX})
$$
\emph{where}
$$
H=\frac{1}{T}\ln \left( \begin{array}{cc} X(T) & 0 \\
0 & Z(T) \end{array} \right),   \eqno (\textrm{X})
$$
\emph{$H \in \mathbf{M}_m(\mathbb{R}), \; Z(T)$ is the monodromy matrix of the
restriction of} (I) \emph {to its periodic invariant manifold  $\,K^d(t),
\;\Phi(t)$ is a periodic matrix with period  $T$ that satisfies the equation }
$$
\frac{d \Phi}{dt}+\Phi H=P(t)\Phi , \eqno (\textrm{XI})
$$
$\Phi(t)\in \mathbf{M}_{n\:m}(\mathbb{R}) $ \emph{for all $t \in\mathbb{R} ,\;
\Phi^+(0)$ is a matrix pseudoinverse to the matrix } $\Phi(0),$ \textit{and} $m=n+d, \;
n\geq d \geq 0.$

Indeed, according to the last remark, the differential equation (I) has the
periodic invariant manifold  $K^d(t)$ on which Eq.~(I) is equivalent to the
second differential equation of system (VIII). Consider the system of
differential equations
\begin{equation}
\label{56}
\frac{dx}{dt}=P(t)x,\qquad   \frac{dz}{dt}=G(t)z,
\end{equation}
which is formed of Eq.~(I) and the second equation of system (VIII).
According to the proof of assertion 3 of Theorem~3, the real canonical form of the monodromy matrix
of this system
\begin{equation}  \label{57}
\left( \begin{array}{cc} X(T) & 0 \\
0 & Z(T) \end{array} \right)
\end{equation}
is the matrix
\begin{equation}  \label{58}
\left( \begin{array}{ccc} A & 0 & 0\\
0 & B & 0 \\ 0 & 0 & B  \end{array} \right),
\end{equation}
where $A$ and $B$ are the blocks of decomposition of the real canonical form of the
matrix  $X(T)$ such that the matrix  $A$ has a real logarithm. Since the matrix
\begin{equation}  \label{59}
\left( \begin{array}{cc} B & 0 \\
0 & B \end{array} \right)
\end{equation}
is formed by pairwise identical Jordan cells, it has a real logarithm. Thus, the
logarithm of matrix (\ref{58}) can be chosen real. Therefore, we can choose the
real logarithm of matrix (\ref{57}) and define the matrix  $H$ according to
relation (X) so that it satisfies the condition  $H \in
\textbf{M}_{n+d}(\mathbb{R}),$ where $d$ is the order of the matrix $B.$ Applying
the Floquet formula (II) to the fundamental matrix of solutions of the system
of differential equations (\ref{56})
\begin{equation}
\left( \begin{array}{cc} X(t) & 0 \\
0 & Z(t) \end{array} \right), \nonumber
\end{equation}
we obtain
\begin{equation}  \label{60}
\left( \begin{array}{cc} X(t) & 0 \\
0 & Z(t) \end{array} \right) =
\left( \begin{array}{cc} \Phi_1(t)  \\
\Phi_2(t)  \end{array} \right)e^{Ht},
\end{equation}
where $H$ is matrix (X) from the space
$\textbf{M}_{n+d}(\mathbb{R}),$ $\Phi_1(t)$ and $\Phi_2(t)$ are periodic
matrices with period  $T,$ and $\Phi_1(t) \in \textbf{M}_{n \: n+d}(\mathbb{R})$
and $\Phi_2(t) \in \textbf{M}_{d \: n+d}(\mathbb{R})$ for all $t\in\mathbb{R}.$
Differentiating equality (\ref{60}), we obtain the following
matrix differential equation for the matrix
$\Phi(t)=\left( \begin{array}{cc} \Phi_1(t)  \\
\Phi_2(t)  \end{array} \right) $:
\begin{equation}
\frac{d\Phi}{dt}+ \Phi H=\left( \begin{array}{cc} P(t) & 0 \\
0 & G(t) \end{array} \right)\Phi. \nonumber
\end{equation}
This equation implies that the matrix $\Phi_1(t)$ satisfies the
differential equation (XI). Finally, multiplying equality (\ref{60}) by the
matrix $\Phi_1^+(0),$ which is pseudoinverse to the matrix $\Phi_1(0),$ we
obtain the equality
\begin{equation}
X(t)=\Phi_1(t)e^{Ht}\Phi_1^+(0), \nonumber
\end{equation}
which coincides (to within notation) with (IX).

\begin{center}\textbf{5. Two Applications of Obtained Results}\end{center}

\textbf{ {1}.} \; Let $x \in \mathbb{R}^n, \:$  let $P(t)$ be a continuous
matrix periodic with period $T,$ let $P(t)\in \textbf{M}_{n}(\mathbb{R})$ for
all $t\in\mathbb{R},\:$ and let $ X(t, x)$ be a function of variables
$t\in\mathbb{R}$ and $x \in \mathbb{R}^n $ that takes values in $\mathbb{R}^n$
and is continuous for all  $t\in\mathbb{R}$ and $x \in \mathbb{R}^n.$

Consider the differential equation \setcounter{equation}{0}
\begin{equation}  \label{1}
\frac{dx}{dt}=P(t)x+ X(t, x).
\end{equation}
Let
\begin{equation}  \label{2}
X(t, x)\equiv 0 .
\end{equation}
Then the differential equation (\ref{1}) has a fundamental matrix of
solutions $X(t)$, which can be represented in the form
\begin{equation}  \label{3}
X(t)=\Phi (t) e^{Ht} \Phi ^+(0),
\end{equation}
and, moreover, the properties of the matrices $\Phi (t)$ and $H$ are determined in
the corollary in the last section.

To simplify the differential equation (\ref{1}), we use
relation (\ref{3}). To this end, we change the variables in (\ref{1}) by
introducing a variable $y \in \mathbb{R}^m$ instead of  $x \in \mathbb{R}^n$
according to the relation
\begin{equation}  \label{4}
x=\Phi (t)y.
\end{equation}
Taking into account that the matrix  $\Phi (t)$ is a solution of the
differential equation (XI), we obtain the following equality from (\ref{1}) and
(\ref{4}):
\begin{equation}  \label{5}
\Phi (t) \left(\frac{dy}{dt}- Hy \right)=X(t, \Phi (t)y).
\end{equation}
We represent this equality in the form
\begin{equation}  \label{6}
\frac{dy}{dt}- Hy= \Phi ^+(t) X(t, \Phi (t)y),
\end{equation}
where $\Phi ^+(t)$ is a matrix pseudoinverse to $\Phi (t)$ that has
the same smoothness and period as  $\Phi (t).$ In particular, as $\Phi ^+(t),$
we can take the first block of the matrix  $(\Phi_1 ^+(t),\: \Phi_2 ^+(t)),$
which is inverse to the matrix $\left( \begin{array}{cc} \Phi_1(t)  \\
\Phi_2(t)  \end{array} \right)$ defined by relation (\ref{60}) of the previous
section. Solving Eq.~(\ref{6}) with respect to $\displaystyle \frac{dy}{dt}- Hy,$ we get
\begin{equation}  \label{7}
\frac{dy}{dt}= Hy + \Phi ^+(t) X(t, \Phi (t)y).
\end{equation}
The selected linear part of Eq.~(\ref{7}) has a constant coefficient matrix,
and the general part preserves the properties of the corresponding part of the
original equation (\ref{1}).

\vspace{6.0mm}
 \textbf{ {2}.} \; Consider the differential equation
\begin{equation}  \label{8}
\frac{dx}{dt}=X(x) + X_1(t, x),
\end{equation}
where $X(x)$ is a continuously differentiable function of  $x$
and  $ X_1(t, x)$ is a continuous function of $t$ and $x $ that
takes values in  $\mathbb{R}^n$ for all $t \in \mathbb{R}$ and $x \in
\mathbb{R}^n, \;  n\geq 2.$

Assume that, under the condition
\begin{equation}  \label{9}
X_1(t, x)\equiv 0,
\end{equation}
Eq.~(\ref{8}) has a $T$-periodic solution
\begin{equation}  \label{10}
x(t)=\xi (\omega t),
\end{equation}
where $\xi(\varphi)$ is a function periodic in $\varphi$ with
period  $2\pi$ and  $\omega =\displaystyle \frac{2\pi}{T}$ is the frequency of
the periodic solution.

The variational equation corresponding to the periodic solution (\ref{10}) of
the differential equation (\ref{8}) with condition (\ref{9}) has the form
\begin{equation}  \label{11}
\frac{d\delta \xi}{dt}= \frac{\partial X(\xi (\omega t))}{\partial x}\,\delta
\xi.
\end{equation}
This equation has the solution
\begin{equation}  \label{12}
\delta \xi=\xi '(\omega t),
\end{equation}
where $'$ stands for the derivative with respect to the variable  $\varphi.$

Indeed, by definition, we have
\begin{equation}  \label{13}
\xi ' (\varphi)\omega=X(\xi (\varphi)).
\end{equation}
Thus,
\begin{equation}  \label{14}
\xi '' (\varphi)\omega=\frac{\partial X(\xi (\varphi
))}{\partial x}\, \xi '(\varphi)
\end{equation}
for all  $\varphi \in \mathbb{R}.$ Substituting $\omega t$ for $\varphi$ in
(\ref{14}), we obtain the identity
\begin{equation}  \label{15}
\frac{d }{dt}\,\xi '(\omega t)=\frac{\partial X(\xi
(\omega t ))}{\partial x}\, \xi '(\omega t),
\end{equation}
which proves the required statement.

Let $B(\varphi)$ be a continuously differentiable periodic matrix with period
$2\pi,$ let $B(\varphi) \in \textbf{M}_{n\,n-1}(\mathbb{R}),$ and let
\begin{equation}  \label{16}
\det (\xi ' (\varphi), B(\varphi))\neq 0
\end{equation}
for all  $\varphi \in \mathbb{R}.$

Using the change of variables
\begin{equation}  \label{17}
\delta \xi= \xi '(\omega t)c+B(\omega t)g,
\end{equation}
we reduce the variational equation (\ref{11}) to the differential equation
\begin{equation}
\xi^{''}(\omega t)\omega c +  \xi^{'}(\omega t)\frac{dc}{dt}+B\:{'}(\omega t)\omega g+
B(\omega t)\frac{dg}{dt}
= \frac{\partial X(\xi(\omega t))}{\partial x}(\xi^{'}(\omega t)c +
B(\omega t)g), \nonumber
\end{equation}
or, with regard for (\ref{14}), to the equation
\begin{equation}
\xi^{'}(\omega t)\frac{dc}{dt}+B(\omega t)\frac{dg}{dt}
=\left(\frac{\partial X(\xi(\omega t))}{\partial x}B(\omega t)-
B\:{'}(\omega t)\omega\right)g. \nonumber
\end{equation}
Solving this equation with respect to the derivatives $\displaystyle\frac{dc}{dt}$ and
$\displaystyle \frac{dg}{dt}$ with the use of the matrix
\begin{equation}  \label{18}
\left( \begin{array}{cc} (\xi ' (\omega t))^+ \\ B^+(\omega t)  \end{array} \right),
\end{equation}
which is inverse to  $(\xi^{'}(\omega t),B(\omega t)),$ we reduce
(\ref{11}) to the system of differential equations
\begin{align} \label{19}
\frac{dc}{dt}=(\xi ' (\omega t))^+
\left(\frac{\partial X(\xi(\omega t))}{\partial x} B(\omega t)-
B\:{'}(\omega t)\omega\right)g, \nonumber \\
\\
\frac{dg}{dt}=B^+(\omega t)\left(\frac{\partial X(\xi(\omega
t))}{\partial x} B(\omega t)- B\:{'}(\omega t)\omega\right)g
\nonumber.
\end{align}

According to the change of variables (\ref{17}), the monodromy matrix of
the system of differential equations (\ref{19}) is similar to the monodromy
matrix of the variational equation (\ref{11}). Thus, the eigenvalues of both
monodromy matrices coincide.

It follows from system (\ref{19}) that one of the eigenvalues of its monodromy
matrix is equal to $1,$ whereas the other eigenvalues are eigenvalues of
the monodromy matrix of the second differential equation of system (\ref{19}). Thus, the same is
true for the eigenvalues of the monodromy matrix of the variational equation
(\ref{11}).

We denote the coefficient matrix of the second differential equation of
system (\ref{19}) by  $Q(\omega t),$ where $Q(\varphi)$ is a periodic matrix
with period $2\pi,$ and consider the differential equation
\begin{equation}\label{20}
\frac{dg}{dt}=Q(\omega t)g.
\end{equation}
By virtue of the corollary in the previous section, the fundamental matrix of
solutions of Eq.~(\ref{20}) $G(t)$ admits the representation
\begin{equation}\label{21}
G(t)=\Phi (\omega t)e^{Ht}\Phi ^{+}(0),
\end{equation}
where
\begin{equation}\label{22}
H=\frac{1}{T} \ln \left( \begin{array}{cc} G(T) & 0 \\ 0 & Z(T)
\end{array} \right) \in \textbf{M}_{m}(\mathbb{R}),
\end{equation}
$2(n-1)\ge m\ge (n-1)\,,\,\,Z(t)$ is the fundamental matrix of the restriction
of the differential equation (\ref{20}) to its periodic invariant manifold
$K^{m-(n-1)}(t),$ $\Phi (\varphi)$ is a periodic matrix with period $2\pi,$
$\Phi (\varphi)\in \textbf{M}_{n-1 \:m}(\mathbb{R})$ for all $\varphi
\in\mathbb{R},$ $\Phi (\varphi)$ satisfies the differential equation
\begin{equation}\label{23}
\frac{d\Phi}{d\varphi}\,\omega +\Phi H=Q(\varphi )\Phi ,
\end{equation}
and $\Phi^{+}(0)$ is a matrix pseudoinverse to $\Phi(0).$

We now use the results obtained above for the introduction of amplitude--phase coordinates
in the neighborhood of the closed curve
\begin{equation}\label{24}
x=\xi(\varphi ),\,\,\varphi \in\mathbb{R} \, ,
\end{equation}
and for the reduction of the differential equation (\ref{8}) in the neighborhood of this
curve to a simpler amplitude--phase system of differential equations.

To this end, we change the variables in Eq.~(\ref{8}) according to the relation
\begin{equation}\label{25}
x=\xi(\varphi)+B(\varphi)g,
\end{equation}
where $B(\varphi)$ is the matrix defined above.

Using equality (\ref{13}), we obtain the following differential
equation instead of (\ref{8}):
\begin{equation}
(\xi^{'}(\varphi )+B^{'}(\varphi )g)\left(\frac{d\varphi}{dt}-\omega\right)
+B(\varphi)\frac{dg}{dt}=\nonumber
\end{equation}
\begin{equation}\label{26}
=X(\xi(\varphi)+B(\varphi)g)-X(\xi(\varphi))-
B^{'}(\varphi)\omega g+
X_{1}(t,\xi(\varphi)+B(\varphi)g).
\end{equation}
We solve Eq.~(\ref{26}) with respect to $\displaystyle\frac{d\varphi}{dt}-\omega$ and
$\displaystyle\frac{dg}{dt}$ by using the matrix
\begin{equation}\label{27}
\left( \begin{array}{cc} L_1(\varphi, g)  \\ L_2(\varphi, g)  \end{array} \right)
\end{equation}
that is inverse to  $(\xi^{'}(\varphi )+B^{'}(\varphi )g,B(\varphi)).$

Choosing a sufficiently small $\delta>0,$ one can easily construct matrix
(\ref{27}) for all
\begin{equation}\label{28}
\varphi\in\mathbb{R},\qquad ||g||\le\delta,
\end{equation}
on the basis of matrix (\ref{18}) by setting
\begin{equation}\label{29}
\left( \begin{array}{cc} L_1(\varphi, 0)  \\ L_2(\varphi, 0)
\end{array} \right)= \left( \begin{array}{cc} (\xi^{'}(\varphi ))^+
\\ B^+(\varphi)  \end{array} \right) .
\end{equation}
Using (\ref{26}), we obtain the system of differential equations
\begin{equation}\label{30}
\frac{d\varphi}{dt}=\omega+L_{1}(\varphi,g)[X(\xi(\varphi)+B(\varphi)g)-X(\xi(\varphi))+
X_{1}(t,\xi(\varphi)+B(\varphi)g)],
\end{equation}
\begin{equation}\label{31}
\frac{dg}{dt}=L_{2}(\varphi,g)[X(\xi(\varphi)+B(\varphi)g)-X(\xi(\varphi))+
X_{1}(t,\xi(\varphi)+B(\varphi)g)].
\end{equation}
We rewrite the differential equation (\ref{31}) in the form
\begin{equation}\label{32}
\frac{dg}{dt}=B^{+}(\varphi)\frac{\partial X(\xi(\varphi))}{\partial x}B(\varphi)g
+G(\varphi,g)+L_{2}(\varphi,g)X_{1}(t,\xi(\varphi)+B(\varphi)g),
\end{equation}
where $G(\varphi,g)$ denotes the function
\begin{equation}
L_{2}(\varphi,g)(X(\xi(\varphi)+B(\varphi)g)-X(\xi(\varphi))-
\frac{\partial X(\xi(\varphi))}{\partial x}B(\varphi)g)+\nonumber
\end{equation}
\begin{equation}\label{33}
+(L_{2}(\varphi,g)-L_{2}(\varphi,0)) \frac{\partial
X(\xi(\varphi))}{\partial x}B(\varphi)g,
\end{equation}
which satisfies the conditions
\begin{equation}\label{34}
G(\varphi,0)=0,\qquad \frac{\partial G(\varphi,0)}{\partial g}=0.
\end{equation}
According to the definition of the matrix $Q(\omega t),$ the coefficient matrix of the selected linear part of the differential equation
(\ref{32}) coincides with the matrix  $Q(\varphi).$ Thus, Eq.~(\ref{32}) takes the form
\begin{equation}\label{35}
\frac{dg}{dt}=Q(\varphi)g+G(\varphi,g)+L_{2}(\varphi,g)X_{1}(t,\xi(\varphi)+B(\varphi)g).
\end{equation}
Let $\Phi(\varphi)$ and $H$ be the matrices determined from representation
(\ref{21}) of a fundamental matrix of solutions of the differential equation
(\ref{20}). With the use of these matrices, we transform the system of
differential equations (\ref{31}), (\ref{32}) by setting
\begin{equation}\label{36}
g=\Phi(\varphi)h.
\end{equation}
As a result, instead of (\ref{35}), we obtain
\begin{equation}
\Phi^{'}(\varphi)h+\Phi(\varphi)\left(\frac{dh}{dt}-Hh\right)+\Phi(\varphi)Hh=
Q(\varphi)\Phi(\varphi)h+G(\varphi,\Phi(\varphi)h)+\nonumber
\end{equation}
\begin{equation}\label{37}
+L_{2}(\varphi,\Phi(\varphi)h)X_{1}(t,\xi(\varphi)+B(\varphi)\Phi(\varphi)h),
\end{equation}
or, taking into account that $\Phi(\varphi)$ is a solution of the differential
equation (\ref{23}),
\begin{equation}\label{38}
\Phi(\varphi)\left(\frac{dh}{dt}-Hh\right)=G(\varphi,\Phi(\varphi)h)
+L_{2}(\varphi,\Phi(\varphi)h)X_{1}(t,\xi(\varphi)+B(\varphi)\Phi(\varphi)h).
\end{equation}
Solving Eq.~(\ref{38}) with respect to $\displaystyle\frac{dh}{dt}-Hh$ with the use of the
matrix $\Phi^{+}(\varphi)$ that is pseudoinverse to $\Phi(\varphi),$ we
obtain
\begin{equation}\label{39}
\frac{dh}{dt}=Hh+\Phi^{+}(\varphi)[G(\varphi,\Phi(\varphi)h)+
L_{2}(\varphi,\Phi(\varphi)h)X_{1}(t,\xi(\varphi)+B(\varphi)\Phi(\varphi)h)].
\end{equation}
By virtue of the results presented above, the change of variables (\ref{36}) reduces the system
of differential equations (\ref{31}), (\ref{32}) to the system
\begin{equation}\label{40}
\frac{d\varphi}{dt}=\omega +f(t,\varphi,\Phi(\varphi)h),
\end{equation}
\begin{equation}\label{41}
\frac{dh}{dt}=Hh+\Phi^{+}(\varphi)F(t,\varphi,\Phi(\varphi)h),
\end{equation}
where $H$ is a matrix of the form (\ref{22}),
\begin{equation}
f(t,\varphi,g)=L_{1}(\varphi,g)[X(\xi(\varphi)+B(\varphi)g)-X(\xi(\varphi))+
X_{1}(t,\xi(\varphi)+B(\varphi)g)],\nonumber
\end{equation}
and
\begin{equation}
F(t,\varphi,g)=L_{2}(\varphi,g)\left[X(\xi(\varphi)+B(\varphi)g)-X(\xi(\varphi))
-\frac{\partial X(\xi(\varphi))}{\partial x}B(\varphi)g\right]+\nonumber
\end{equation}
\begin{equation}
+(L_{2}(\varphi,g)-L_{2}(\varphi,0))\frac{\partial X(\xi(\varphi))}{\partial x}B(\varphi)g.\nonumber
\end{equation}
The system of differential equations (\ref{40}), (\ref{41}) is the required one.

Thus, by using the superposition of changes (\ref{25}) and (\ref{36}), and,
hence, the change of variables
\begin{equation}
x=\xi(\varphi)+B(\varphi)\Phi(\varphi)h,\nonumber
\end{equation}
the original differential equation (\ref{8}) can be reduced in the neighborhood of the closed
curve (\ref{24}) to the system of differential equations (\ref{40}),
(\ref{41}), where the functions  $f(t,\,\varphi,\,g)$ and  $F(t,\,
\varphi,\,g)$ are continuous in the variables  $t,$ $\varphi,$  and  $g$ for
$t\in\mathbb{R},$ $\varphi\in\mathbb{R},$ and $g\in\mathbb{R}^{n-1},$ $||g||\le\delta,$ take values in $\mathbb{R}$ and
$\mathbb{R}^{n-1},$ respectively, and are periodic in $\varphi$ with
period  $2\pi,$ the matrix  $\Phi(\varphi)$ belongs to $\mathbf{M}_{n-1 \:m}(\mathbb{R})$
for all  $\varphi\in\mathbb{R}$ and is periodic with period  $2\pi,$ the matrix
$H$ belongs to $\mathbf{M}_{m}(\mathbb{R}) ,$ its eigenvalues are the numbers
\begin{equation}
\frac{1}{T}\ln \lambda_{j},\qquad j=\overline{1, \,n-1},\nonumber
\end{equation}
and their  $p$-fold repetitions, $1\ge p_{j}\ge
0,\,\,\sum\limits_{j=1}^{n-1}p_{j}=m-(n-1),\,\, 2(n-1)\ge m\ge (n-1),$ and $1$ and
$\lambda_{1},\ldots,\lambda_{n-1}$ are the eigenvalues of the monodromy matrix of
the variational equation (\ref{11}).

The reduction of the differential equations considered above to equations with
the constant matrix of coefficients in their separated linear part is essential for
the subsequent investigation of these equations. A confirmation of this statement can be found, e.g., in
[\,5,\,6\,], where, however, the problem of this reduction was only
partially solved.

\vspace*{36pt}

\noindent 1. A. M. Samoilenko, \textit{On Invariant Manifolds of Linear
Differential Equations} [in Ukrainian], Preprint No.~2009.7, Institute of
Mathematics, Ukrainian National Academy of Sciences, Kyiv (2009). \\
2. G. Floquet, ``Sur les \'{e}quations diff\'{e}rentielles lin\'{e}aires \`{a}
coefficients p\'{e}riodiques,'' \textit{ Ann. \'{E}cole Norm. Sup.,} No.~12,
47--88 (1883). \\
3. F. R. Gantmakher, \textit{Theory of Matrices} [in Russian], Nauka, Moscow (1988). \\
4. A. M. Samoilenko, \textit{Elements of the Mathematical Theory of
Multifrequency Oscillations.  Invariant Tori} [in Russian], Nauka, Moscow
(1987). \\
5.  N. N. Bogolyubov and Yu. A. Mitropol'skii, \textit{Asymptotic Methods in
the Theory of Nonlinear Oscillations} [in Russian], Fizmatgiz, Moscow (1963);
N. N. Bogolyubov, \textit{Collection of Scientific Works} [in Russian], Vol.~3,
Nauka, Moscow (2005) \\
6. A. M. Samoilenko and L. Recke, ``Conditions for synchronization
 of one oscillating system,'' \textit{Ukr. Mat. Zh.,} \textbf{57},
No.~7, 922 -- 945 (2005). \\

\end{document}


\newcounter{abcd}
\setcounter{abcd}{1}

\fontsize{14.4}{18pt}\selectfont  

УДК 517.9 \vspace{6.0mm}

 \textbf{А. М. Самойленко}\quad
(Ін-т математики НАН України, Київ)\vspace{6.0mm}

\textbf{Про інваріантні многовиди лінійних диференціальних
рівнянь.~II} \vspace{8.0mm}

Продовжимо дослідження, розпочаті в роботі [\,1\,]. \vspace{6.0mm}

\begin{center}\textbf{3. Еквівалентність лінійних диференціальних
рівнянь різних порядків.}\end{center}

Як і в [\,1\,], позначимо через $\Phi_{1}(t)$ та $\Phi_{2}(t)$ блоки
невиродженої неперервно диференційовної для всіх $t\in \mathbb{R}$
матриці
\begin{equation}
\Phi(t)=\left(\begin{array}{cc}\Phi_{1}(t)\\\Phi_{2}(t)\end{array}\right),
\nonumber \end{equation} через $\Phi_{1}^+ (t)$ та $\Phi_{2}^+ (t)$
--- блоки оберненої до $\Phi(t)$ матриці
\begin{equation}
\Phi^{-1}(t)=(\Phi_{1}^+ (t), \Phi_{2}^+ (t)), \nonumber
\end{equation}
\begin{equation}
\Phi_{1}(t)\in \textbf{M}_{n \,m}(\mathbb{R}),\qquad \Phi_{2}(t)\in
\textbf{M}_{m-n \,m}(\mathbb{R}), \nonumber \end{equation}
\begin{equation}
\Phi_{1}^+(t)\in \textbf{M}_{m \,n}(\mathbb{R}),\qquad
\Phi_{2}^+(t)\in \textbf{M}_{m \,m-n}(\mathbb{R}), \nonumber
\end{equation}
 через
$M_{1}(t)=\Phi_{1}^{+}(t)\Phi_{1}(t)$ та
$M_{2}(t)=\Phi_{2}^{+}(t)\Phi_{2}(t)$ --- проектори рангів $n$ та
$m-n$ відповідно, через $M^n(t)$ та $M^{m-n}(t)$ --- гіперплощини
\begin{equation}
 M^{n}(t)=\{y\in \mathbb{R}^{m}: y=M_{1}(t)y\},
\nonumber \end{equation}
\begin{equation}
 M^{m-n}(t)=\{y\in \mathbb{R}^{m-n}: M_{1}(t)y=0\}
\nonumber \end{equation} розмірів $n$ та $m-n$ відповідно, через
$L(M, Q)$
--- матричний оператор вигляду
\begin{equation}
L(M, Q)=\frac{dM}{dt} +MQ-QM. \nonumber \end{equation}

\textbf{Теорема 2.} {\em Якщо підпростори \/} $M^{n}(t) $ {\em та
\/} $M^{m-n}(t), $ {\em взяті разом, є ін\-ва\-рі\-ант\-ними многовидами
диференціального рівняння \/}
\begin{equation}\label{rim1}
\frac{dy}{dt}=Q(t)y, \tag{\Roman{abcd}}
\end{equation}
\addtocounter{abcd}{1}
{\em то заміна змінних  \/}
\begin{equation}\label{rim2}
y=\Phi_{1}^+(t)x+\Phi_{2}^+(t)z \tag{\Roman{abcd}}
\end{equation}
\addtocounter{abcd}{1}
{\em зводить це рівняння до системи
диференціальних рівнянь \/}
\begin{equation}\label{rim3}
\frac{dx}{dt}=P(t)x, \qquad \frac{dz}{dt}=G(t)z \tag{\Roman{abcd}}
\end{equation}
\addtocounter{abcd}{1} {\em з матрицями коефіцієнтів  \/}
\begin{equation}\label{rim4}
P(t)= \left (\frac{d\Phi_1(t)}{dt}+\Phi_1(t)Q(t)
\right)\Phi_1^{+}(t), \tag{\Roman{abcd}}
\end{equation}
\addtocounter{abcd}{1}
\begin{equation}\label{rim5}
G(t)= \left (\frac{d\Phi_2(t)}{dt}+\Phi_2(t)Q(t)
\right)\Phi_2^{+}(t), \tag{\Roman{abcd}}
\end{equation}
\addtocounter{abcd}{1} {\em і навпаки: якщо диференціальне рівняння}
(\textrm{I}) {\em заміною змінних} (\textrm{II}) {\em зводиться  до
системи диференціальних рівнянь}
(\textrm{III}), {\em то підпростори \/} $M^{n}(t) $ {\em та \/} \\
$M^{m-n}(t), $ {\em взяті разом, є інваріантними многовидами
рівняння} (\textrm{I}), {\em а матриці коефі\-ці\-єн\-тів системи}
(\textrm{III}) {\em визначаються формулами}  (\textrm{IV}),
(\textrm{V}).

Дійсно, нехай підпростори  $M^{n}(t) $  та  $M^{m-n}(t), $  взяті
разом, є інваріантними многовидами рівняння (I). Тоді згідно з
твердженнями 2 основної теореми роботи [\,1\,] на $M^{n}(t) $ та
$M^{m-n}(t) $ рівняння (I) еквівалентне відповідно першому та
другому з рівнянь системи (III), причому з матрицями коефіцієнтів
цих рівнянь, визначених формулами (IV) та (V) відповідно. Нехай
$Y(t), X(t) $  та  $Z(t),$ де $Y(0)=E,\, X(0)=E$ та $Z(0)=E \,,$ ---
фундаментальні матриці розв'язків рівнянь (I) та (III), $E$ ---
одиничні матриці відповідних розмірів. Згідно з твердженнями 2
основної теореми роботи [\,1\,] маємо співвідношення
\begin{equation}
\label{1} Y(t)\Phi^{+}_1(0)=\Phi^{+}_1(t)X(t), \qquad
Y(t)\Phi^{+}_2(0)=\Phi^{+}_2(t)Z(t)
\end{equation}
для всіх $t\in \mathbb{R}.$ Отже, у відповідності до (\ref{1})
\begin{equation}
\label{2} Y(t)(\Phi^{+}_1(0), \Phi^{+}_2(0))=
(\Phi^{+}_1(t),\Phi^{+}_2(t))\left(\begin{array}{cc} X(t)&0
\\ 0 & Z(t)\end{array}\right)
\end{equation}
для всіх $t\in \mathbb{R}.$ Із (\ref{2}) випливає рівність
\begin{equation}
\label{3}Y(t)=\Phi^{+}_1(t)X(t)\Phi_1(0) +
\Phi^{+}_2(t)Z(t)\Phi_2(0)
\end{equation}
для всіх $t\in \mathbb{R}.$ Таким чином, для довільного  $y_0 \in
\mathbb{R}^m$
\begin{equation}
\label{4}Y(t)y_0=\Phi^{+}_1(t)X(t)x_0 + \Phi^{+}_2(t)Z(t)z_0
\end{equation}
для всіх $t\in \mathbb{R}$ та $x_0, z_0$, вибраних згідно з умовою
\begin{equation} x_0=\Phi_1(0)y_0, \qquad z_0=\Phi_2(0)y_0.
\nonumber \end{equation} Рівність (\ref{4}) якраз і означає, що
заміна змінних (II) зводить диференціальне рівняння (I) до системи
диференціальних рівнянь (III).

Нехай тепер диференціальне рівняння (I) заміною змінних (II)
зводиться до системи диференціальних рівнянь (III). Оскільки
підпростори  $z=0$ та $x=0$ є інваріантними многовидами системи
(III), то з цього та заміни (II) випливають співвідношення
(\ref{1}), а з них --- співвідношення
\begin{equation}
\label{5} X(t)=\Phi_1(t)Y(t)\Phi_1^{+}(0), \qquad
Z(t)=\Phi_2(t)Y(t)\Phi_2^{+}(0)
\end{equation}
для всіх $t\in \mathbb{R}.$

Підставляючи (\ref{5}) у формули (\ref{1}),  отримаємо, що
\begin{equation}
\label{6} Y(t)\Phi^{+}_1(0)=M_1(t)Y(t)\Phi^{+}_1(0), \qquad
Y(t)\Phi^{+}_2(0)=M_2(t)Y(t)\Phi^{+}_2(0)
\end{equation}
для всіх $t\in \mathbb{R}.$

Із першого співвідношення (\ref{6}) випливає, що
\begin{equation}
\label{7} y(t)=M_1(t)y(t)
\end{equation}
для будь-якого розв'язку $y(t)$ рівняння (I),  який задовольняє умову
\begin{equation}
\label{8} y(0)=\Phi^{+}_1(0)c,
\end{equation}
де $c$ --- довільна стала з $\mathbb{R}^n.$ Оскільки точки (\ref{8})
заповнюють підпростір $M^n(0)$, то згідно з (\ref{7}) інтегральні
криві $(t, y(t))$ рівняння (I), що при $t=0$ проходять через точки
підпростору $M^n(0)$, належать підпростору $M^n(t)$ для будь-якого
$t\in \mathbb{R}.$ Цього достатньо, щоб підпростір $M^n(t)$ був
інваріантним многовидом рівняння (I).

Із другого співвідношення (\ref{6}) випливає, що
\begin{equation}
\label{9} y(t)=M_2(t)y(t)
\end{equation}
для будь-якого розв'язку $y(t)$ рівняння (I), який задовольняє умову
\begin{equation}
\label{10} y(0)=\Phi^{+}_2(0)c,
\end{equation}
де $c$ --- довільна стала з $\mathbb{R}^{m-n}.$

Аналогічні до наведених вище міркування доводять, що підпростір
\begin{equation}
M_2^{m-n}(t)=\{y \in \mathbb{R}^{m-n}: \; y=M_2(t)y \} \nonumber
\end{equation} є інваріантним многовидом рівняння (I).

Згідно з лемою 1 роботи [\,1\,] має місце рівність
\begin{equation} M_2^{m-n}(t)= M^{m-n}(t) \nonumber \end{equation}
для будь-якого $t\in \mathbb{R}.$ Це доводить, що підпростір
$M^{m-n}(t)$ є інваріантним многовидом диференціального рівняння
(I). Таким чином, підпростори $M^{n}(t)$ та $M^{m-n}(t)$, взяті
разом, є інваріантними многовидами рівняння (I). Згідно з
твердженнями 2 основної теореми [\,1\,] цього достатньо для того,
щоб справджувалися формули (IV) та (V).

Нехай $F(t) \in M_{p \;n}(\mathbb{R}),\, n>p,\, F^+(t) \in M_{p
\;n}(\mathbb{R}),\, \textrm{rank} F(t)=p, \,F(t)$ та $F^+(t)$
 є неперервно диференційовними функціями для всіх $t\in
\mathbb{R}.$ Крім того, нехай $F^+(t)$ є матрицею, псевдооберненою
для матриці $F(t), K(t)=F^+(t)F(t).$ Нарешті, нехай підпростір
\begin{equation}
K^p(t)=\{ x \in \mathbb{R}^{n}: \; x=K(t)x \} \nonumber
\end{equation} є інваріантним многовидом диференціального рівняння
\begin{equation}\label{rim6}
\frac{dx}{dt}=P(t)x, \tag{\Roman{abcd}}
\end{equation}
\addtocounter{abcd}{1}
і це рівняння на $K^p(t)$ еквівалентне
диференціальному рівнянню
\begin{equation}\label{rim7}
\frac{dz}{dt}=R(t)z. \tag{\Roman{abcd}}
\end{equation}
\addtocounter{abcd}{1}

Назвемо:

систему диференціальних рівнянь (III) \textit{декомпозицією}
диференціального рівняння (I), якщо заміна змінних (II) зводить це
рівняння до системи  рівнянь (III);

диференціальне рівняння (VII) \textit{звуженням} диференціального
рівняння (VI) на підпростір $K^p(t)$, якщо підпростір $K^p(t)$ є
інваріантним многовидом рівняння (VI) і на $K^p(t)$ це рівняння
еквівалентне диференціальному рівнянню (VII).

Диференціальні рівняння (I) та (VI) будемо називати еквівалентними,
якщо рівнян\-ня (VI) разом з його звуженням на $K^{m-n}(t)$ (VII) є
декомпозицією рівняння (I).

За визначенням фундаментальні матриці розв'язків еквівалентних
дифе\-рен\-ці\-аль\-них рівнянь взаємно визначаються одна через одну
та через матриці, які визначають інваріантні підпростори цих
диференціальних рівнянь. Дійсно, згідно з наведеними вище означеннями та
з урахуванням того, що
\begin{equation}\label{11} G(t)=R(t)
\end{equation}
для всіх $t\in \mathbb{R},$ справедливими  є співвідношення
(\ref{1}) та
\begin{equation}\label{12} X(t)F^+(0)=F^+(t)Z(t)
\end{equation}
між фундаментальними матрицями розв'язків $Y(t),$ $X(t)$ та $Z(t)$
диферен\-ці\-аль\-них рівнянь (I), (VI) та (VII).

Із (\ref{1}) та (\ref{12}) випливає, що
\begin{equation}
Y(t)=(\Phi_{1}^{+}(t)X(t)+\Phi_{2}^{+}(t)F(t)X(t)F^{+}(0))
\left(\begin{array}{cc}\Phi_{1}(0)\\\Phi_{2}(0)\end{array}\right)=\nonumber
\end{equation}
\begin{equation}\label{13}
=\Phi_{1}^{+}(t)X(t)\Phi_{1}(0)+\Phi_{2}^{+}(t)F(t)X(t)F^{+}(0)\Phi_{2}(0),
\end{equation}
\begin{equation}\label{14}
X(t)=\Phi_{1}(t)Y(t)\Phi_{1}^{+}(0)
\end{equation}
для всіх $t\in \mathbb{R} .$ Співвідношення (\ref{13}), (\ref{14})
визначають взаємозв'язок між фундаментальними матрицями розв'язків
еквівалентних диференціальних рівнянь (I) та (VI).

Поняття еквівалентності диференціальних рівнянь $m$-го та $n$-го
порядків ви\-зна\-че\-но вище для
\begin{equation} m>n>m-n. \nonumber \end{equation}
Його легко поширити на випадок, коли
\begin{equation}\label{15}
m=2n.
\end{equation}
Дійсно, оскільки простір $\mathbb{R}^{n}$ є інваріантним многовидом
диференціального рівняння (VI) і на ньому рівняння (VI) еквівалентне
диференціальному рів\-нян\-ню (VIІ) з тією самою матрицею коефіцієнтів,
то у випадку (\ref{15}) еквіва\-лент\-ність диференціальних рівнянь (І)
та (VI) визначається декомпозицією рівняння (І) до системи рівнянь
\begin{equation} \frac{dx}{dt}=P(t)x,\qquad \frac{dz}{dt}=P(t)z. \nonumber \end{equation}

Із наведених вище результатів випливає наступне твердження.

\textbf{Наслідок.} \emph{Диференціальні рівняння} (І) \emph{та} (VI)
\emph{еквівалентні тоді й лише тоді, коли}
$$ L(M(t), Q(t))=0, \qquad L(K(t), P(t))K(t)=0,  \eqno (\textrm{VIII}) $$
$$ P(t)=\left (\frac{d \Phi_1(t)}{dt}+\Phi_1(t)Q(t) \right )\Phi_1^+(t),  \eqno (\textrm{XI}) $$
$$ \left (\frac{d \Phi_2(t)}{dt}+\Phi_2(t)Q(t) \right )\Phi_2^+(t)=
 \left (\frac{d F(t)}{dt}+F(t)P(t) \right )F^+(t)  \eqno (\textrm{X}) $$
\emph{для всіх} $t \in\mathbb{R}$.

 Справді, нехай диференціальні рівняння (І) та (VI) еквівалентні.
Тоді має місце декомпозиція рівняння (І) до системи (ІІІ), друге з
рівнянь якої є звуженням на $K^{m-n}(t)$ диференціального рівняння
(VI). Із визначення де\-ком\-по\-зи\-ції та теореми 2 випливає, що
підпростори $M^{n}(t)$ та $M^{m-n}(t)$ є інварі\-ант\-ни\-ми
многовидами диференціального рівняння (І), а із визначення звуження
диференціального рівняння (VI) на підпростір $K^{m-n}(t)$ випливає,
що $K^{m-n}(t)$ є інваріантним многовидом цього рівняння. Згідно з
твердженнями 1 та 2 основної теореми роботи [\,1\,] цього достатньо,
щоб виконувались співвідношення (VIII) та (IX). Крім того, цього
достатньо, щоб матриці коефіцієнтів  диферен\-ці\-аль\-них рівнянь
(І), (ІІІ) та (VII) задовольняли співвідношення
\begin{equation}
\label{16}
 G(t)= \left (\frac{d \Phi_2(t)}{dt}+ \Phi_2(t)Q(t) \right )
 \Phi_2^+(t),
\end{equation}
\begin{equation}
\label{17}
R(t)=\left (\frac{d F(t)}{dt}+F(t)P(t) \right )F^+(t)
\end{equation}
та
\begin{equation}
\label{18}  G(t)=R(t)
\end{equation}
 для всіх $t \in\mathbb{R}$.

Останнє з наведених вище співвідношень доводить рівність (Х).

Нехай виконуються співвідношення (VIII) -- (Х). Тоді згідно з
тверджен\-ня\-ми 1 та 2 основної теореми роботи [\,1\,] підпростори
$M^{n}(t)$ та $M^{m-n}(t)$ є інварі\-ант\-ни\-ми многовидами
диференціального рівняння (І), а підпростір $K^{m-n}(t)$ ---
інваріантним многовидом диференціального рівняння (VI), при цьому
матриці коефіцієнтів відповідних диференціальних рівнянь $G(t)$ та
$R(t)$ визначаються співвідношеннями (\ref{16}) та (\ref{17}), отже,
згідно з умовою (X) задовольняють рівність (\ref{18}). Звідси
випливає, що у відповідності з теоремою 2 декомпозицією
диференціального рівняння (I) є система диференціальних рівнянь
(III), друге з рівнянь якої є  звуженням диференціального рівняння
(VI) на підпростір $K^{m-n}(t).$ Це доводить, що із співвідношень
(VIII) -- (Х) випливає еквівалентність диференціальних рівнянь (І)
та (VI) згідно з визна\-чен\-ням останньої.

Відзначимо, що при $m=2n$ умови (VIII) -- (Х) спрощуються за рахунок
того, що в даному випадку матриці  $F(t)$ та $K(t)$ є одиничними, та
зводяться до співвідношень
  вигляду
\begin{equation}
 L(M(t), Q(t))=0, \qquad
  P(t)=\left (\frac{d \Phi_1(t)}{dt}+\Phi_1(t)Q(t)
\right )\Phi_1^+(t) = \nonumber \end{equation}
\begin{equation}= \left (\frac{d
\Phi_2(t)}{dt}+\Phi_2(t)Q(t) \right )\Phi_2^+(t)
  \nonumber \end{equation}
для будь-якого $t \in\mathbb{R}$.

Зауважимо також, що еквівалентність диференціальних рівнянь (І) та
(VI) означає справедливість співвідношень
\begin{equation}
\label{19}
Y(t) \Phi_{1}^{+}(0)=\Phi_{1}^{+}(t) X(t), \qquad Y(t)
\Phi_{2}^{+}(0)= \Phi_{2}^{+}(t) F(t) X(t) F^{+}(0)
\end{equation}
між фундаментальними матрицями розв'язків рівнянь (І) та (VI) $Y(t)$
та $X(t)$, а також усіх інших співвідношень, які можна отримати з
(\ref{19}) перетвореннями останніх.

\begin{center}\textbf{4. Додаток до теорії Флоке -- Ляпунова.}\end{center}

Розглянемо лінійне диференціальне рівняння
 \setcounter{abcd}{1}
\begin{equation}\label{rim1}
\frac{dx}{dt}=P(t)x, \tag{\Roman{abcd}}
\end{equation}
\addtocounter{abcd}{1}
де $ x \in\mathbb{R}^n, \, P(t)\in
\textbf{M}_n(\mathbb{R}), \, P(t)$
--- неперервна та періодична з періодом $T$ матриця.

Фундаментальна матриця розв'язків рівняння (І) $X(t), \, X(0)=E, $
згідно з відомою теоремою Г. Флоке [\:2\,] має представлення
$$ X(t)=\Phi (t) e^{Ht},  \eqno(\textrm{II})$$
де $\Phi(t)$ --- періодична за змінною $t$ з періодом $T$ матриця,
$H$ --- стала матриця, що визначається  матрицею монодромії $X(T)$
рівняння (І) за формулою
$$ H=\frac{1}{T}\ln X(T).  \eqno(\textrm{III})$$
Логарифм --- багатозначна функція, і його дійсне значення не завжди
існує. Отже, формула (І) з матрицею (III) такою, що
$$ H \in \textbf{M}_n(\mathbb{R}),  \eqno(\textrm{ІV})$$
не завжди справджується. Згідно з теорією матриць [\,3\,] умова (ІV)
виконується тоді й лише тоді, коли кожен елементарний дільник, що
відповідає від'ємним власним значенням матриці $X(T),$ повторюється
парну кількість разів. Отже, лише у цьому випадку рівність (II)
справджується з матрицями $\Phi(t)$ та $H$ із простору дійсних
матриць $\textbf{M}_n(\mathbb{R}).$

У випадку, коли умову (ІV) не можна задовольнити, представлення Г.
Флоке (II) існує лише з матрицями $\Phi(t)$ та $H$  із простору
$\textbf{M}_n(\mathbb{C}),$ де $\mathbb{C}$ --- площина комплексних
чисел, або ж це представлення трансформується у рівність (II) з
дійсними матрицями $\Phi(t)$ та $H,$ перша з яких є періодичною з
періодом $2T,$ а друга визначається формулою
$$ H=\frac{1}{2T}\ln X(2T).  \eqno(\textrm{V})$$
Представлення Г. Флоке (II) з матрицею (V) --- це наслідок наявності
серед власних значень матриці монодромії рівняння (I) від'ємних
чисел.

Розглянемо детальніше диференціальне рівняння (I)  саме із такою
власти\-віс\-тю його матриці монодромії та доведемо для нього низку
тверджень, неві\-до\-мих раніше.

\textbf{Теорема 3.} {\em Нехай матриця коефіцієнтів диференціального
рівняння} (I)  $ P(t)\in \textbf{M}_n(\mathbb{R})$ {\em для
будь-якого $t \in\mathbb{R} $, неперервна на  $\mathbb{R}$ та
періодична за $t $ з періодом $T.$ \/}

 {\em Тоді:\/}

1. {\em Алгебраїчна кількість від'ємних власних значень матриці
монодромії $X(T)$ рівняння} (I) $\,p$ {\em є парним числом. \/}

2. {\em Рівність} (II) {\em справджується для матриці \/}
$$ H=\frac{1}{T}\ln (X(T)I),  \eqno(\textrm{VI})$$
{\em  де $I$ --- дійсна матриця, визначена умовами  \/}
$$ I^2=E, \qquad \ln (X(T)I) \in \textbf{M}_n(\mathbb{R}),  $$
{\em та періодичної матриці  $\Phi(t)$ такої, що \/}
$$ \Phi(t+T)I_1=\Phi(t)I_1, \qquad   \Phi(t+T)I_2=-\Phi(t)I_2  \eqno(\textrm{VII})$$
{\em для всіх $t \in\mathbb{R},$  де \/}
$$ I_1=\frac{E+I}{2}, \qquad I_2=\frac{E-I}{2}.  $$

3. {\em Існує невироджена неперервно диференційовна та дійсна для
всіх $t \in\mathbb{R}$ матриця $(U(t), \, V(t))$, періодична  з
періодом $T$ і така, що заміна змінних
$$ x= U(t)z_1+ V(t)z_2  $$
зводить диференціальне рівняння} (I) {\em до системи диференціальних
рівнянь \/}
$$\frac{d z_1}{dt}=H_1 z_1, \qquad \frac{d z_2}{dt}=G(t) z_2,  \eqno(\textrm{VIII})$$
{\em у якої $H_1$ --- стала, $G(t)$ --- періодична з періодом $T$
матриці та множина власних значень матриці монодромії $Z_2(T)$
другого її рівняння  є або  мно\-жи\-ною всіх від'ємних власних
значень матриці $X(T)$, або  її підмножиною. \/}

Для доведення теореми використаємо представлення матриці $X(T)$
 її жорда\-но\-вою формою $J(\lambda)$
\begin{equation}
 X(T)=S J(\lambda)S^{-1} \nonumber \end{equation}
та отримаємо рівність \setcounter{equation}{0}
\begin{equation}
\label{1} \det X(T)=\prod ^n_{\nu =1} \lambda _ \nu,
\end{equation}
що пов'язує визначник матриці $X(T)$ з її власними значеннями
$\lambda _ \nu\/, \, \nu=\overline{1,\, n}.$

Тепер застосуємо  формулу Ліувіля -- Остроградського -- Якобі та
подамо визначник матриці $X(T)$ через слід матриці коефіцієнтів
рівняння (I)
\begin{equation}
\label{2}
 \det X(T)= \exp \left \{ \textrm{tr} \,P(t) dt \right \} .
\end{equation}
Прирівнюючи праві частини формул (\ref{1}) та (\ref{2}), отримаємо
рівність, яка доводить, що
\begin{equation}
\label{3}\prod ^n_{\nu =1} \lambda _ \nu > 0.
\end{equation}
Оскільки з кожної пари комплексно-спряжених власних значень матриці
$X(T)$ у добутку всіх її власних значень отримується додатне число,
то з (\ref{3}) випливає, що добуток всіх від'ємних власних значень
матриці $X(T)$ також є додатним числом. Отже, алгебраїчна кількість
від'ємних власних значень матриці $X(T),$ тобто сума, яку складають
кратності як корені характеристичного рівняння
 всіх різних від'ємних власних значень матриці
$X(T),$  є парним числом.

Перш ніж доводити твердження 2 теореми  3, зауважимо, що коли
матриця $X(T)$ має дійсний логарифм, то, поклавши $I=E$, можна
звести рівності (II), (VI) до формул Г.~Флоке (II), (III) з
 матрицею $\Phi(T),$ яка має властивості, що випливають
з цих формул та вказані у твердженнях 2 теореми 3.

Отже, залишилось розглянути випадок, коли матриця $X(T)$ має
від'ємні власні значення та не має дійсного логарифма. В цьому
випадку дійсну канонічну форму матриці $X(T)$ можна подати у вигляді
розбиття на два блоки $A$ та $B$, причому $A$ або порожній, або має
дійсний логарифм, а $B$ має лише від'ємні власні значення та не має
дійсного логарифма.

Нехай $B\in \textbf{M}_d(\mathbb{R}),$ де
\begin{equation}
\label{4}n>d.
\end{equation}
Тоді виконується рівність
\begin{equation}
\label{5} X(T)=S \left(\begin{array}{cc} A&0
\\ 0 & B\end{array}\right) S^{-1}\, ,
\end{equation}
де $S,\,\,A$ та $B$ --- дійсні матриці з вказаними вище
властивостями для $A$ та $B$.

Покладемо
\begin{equation}
\label{6}Y(t)=S^{-1}X(t)S,\qquad B_{1}=-B.
\end{equation}
Оскільки згідно з властивостями фундаментальної матриці розв'язків \\
рівняння (I)
\begin{equation}
\label{7}X(t+T)=X(t)X(T)\,,
\end{equation}
то із (\ref{5}), (\ref{6}) та (\ref{7}) випливає, що
\begin{equation}
\label{8}Y(t+kT)=S^{-1}X(t)X^{k}(T)S=S^{-1}X(t)SS^{-1}X^{k}(T)S=\nonumber
\end{equation}
\begin{equation} =Y(t) \left(\begin{array}{cc} A^k&0
\\ 0 & B^k\end{array}\right)=Y(t)\left(\begin{array}{cc} A^k&0
\\ 0 & (-1)^k B_1^k\end{array}\right)
\end{equation}
для будь-якого цілого $k$.

Подамо $Y(t)$ у блочному вигляді
\begin{equation}\label{9}
Y(t)=(Y_{1}(t),Y_{2}(t)),
\end{equation}
узгодженому з розбиттям (5) матриці $X(T)$ на блоки $A$ та $B$ .
Тоді  з формул (\ref{8}) одержимо співвідношення
\begin{equation}\label{10}
Y_{1}(t+kT)=Y_{1}(t)A^{k},\qquad
Y_{2}(t+kT)=(-1)^{k}Y_{2}(t)B_{1}^{k}
\end{equation}
для будь-якого цілого $k$.

Оскільки власні числа матриці $B_{1}$ в силу визначення цієї матриці
(\ref{6}) додатні, то обидві матриці $A$ та $B_{1}$ мають дійсні
логарифми $\ln A$ та $\ln B_{1}$.

Враховуючи  викладене вище, із (\ref{10}) отримуємо, що
\begin{equation}
Y_{1}(t)=Y_{1}\left(t-\left[\frac{t}{T}\right]T+\left[\frac{t}{T}\right]T\right)=Y_{1}\left
(t- \left[\frac{t}{T}\right]T \right)A^{\left[
\frac{t}{T}\right]}=\nonumber
\end{equation}
\begin{equation}\label{11}
=Y_{1}\left(t- \left[\frac{t}{T}\right]T \right)\exp \left \{
\left(\left[ \frac{t}{T}\right]T-t\right) \frac{\ln A}{T} \right \}
 \exp \left \{ \frac{t}{T} \ln A \right \}  ,
\end{equation}
\begin{equation}
Y_{2}(t)=Y_{2}\left(t-\left[\frac{t}{T}\right]T+\left[\frac{t}{T}\right]T\right)=Y_{2}\left(t-
\left[\frac{t}{T}\right]T\right)(-1)^{\left[\frac{t}{T}\right]}B_{1}^{\left[\frac{t}{T}\right]}=\nonumber
\end{equation}
\begin{equation}\label{12}
=(-1)^{\left[\frac{t}{T}\right]}Y_{2}\left(t-
\left[\frac{t}{T}\right]T\right) \exp \left \{
\left(\left[\frac{t}{T}\right]T-t\right)\frac{\ln B_{1}}{T} \right
\} \exp \left \{\frac{t}{T}\ln B_{1} \right \}
\end{equation}
для всіх $t\in\mathbb{R} ,$ де $[\,t\,]$ означає цілу частину числа
$t$.

Позначимо множники перед $\, \exp \left \{\displaystyle
\frac{t}{T}\ln A \right \}$ та $\, \exp \left \{\displaystyle
\frac{t}{T}\ln B_{1} \right \}$ у формулах (\ref{11}) та (\ref{12})
через $\Phi_{1}(t)$ та $\Phi_{2}(t)$ відповідно. Тоді згідно з
(\ref{9}), (\ref{11}) та (\ref{12}) одержимо рівність
\begin{equation}\label{13}
Y(t)=(\Phi_{1}(t),\Phi_{2}(t)) \left(\begin{array}{cc} \exp \left
\{\displaystyle \frac{t}{T}\ln A \right \} &0
\\ 0 & \exp \left \{\displaystyle \frac{t}{T}\ln B_{1} \right \} \end{array}\right)
\end{equation}
для всіх $t\in\mathbb{R}$, з якої випливає, що матриці $\Phi_{1}(t)$
та $\Phi_{2}(t)$ неперервно диференційовні на $\mathbb{R}$. Крім
того, як видно із введених позначень,  матриця $\Phi_{1}(t)$ є
періодичною  з періодом $T$, а матриця $\Phi_{2}(t)$ як добуток
функції $(-1)^{[ \frac{t}{T}]}$ на періодичну з періодом $T$ матрицю
задовольняє умову
\begin{equation}\label{14}
\Phi_{2}(t+T)=-\Phi_{2}(t)
\end{equation}
для всіх $t\in\mathbb{R}$.

Позначимо через $I_{0}$ матрицю
$$\left( \begin{array}{cc} E_1 & 0 \\ 0 & -E_2 \end{array} \right),$$
де $E_{1}$ та $E_{2}$ --- одиничні матриці із
$\textbf{M}_{n-d}(\mathbb{R})$ та $\textbf{M}_{d}(\mathbb{R})$
відповідно. Тоді
$$Y(T)I_{0}=\left( \begin{array}{cc} A & 0 \\ 0 & B_1 \end{array} \right),$$
і формула (\ref{13}) набуває вигляду
\begin{equation}\label{15}
Y(t)=(\Phi_{1}(t),\Phi_{2}(t)) \exp \left \{ \frac{t}{T}\ln (Y(T)I_{0}) \right \}.
\end{equation}
Із (\ref{15}), враховуючи першу з рівностей  (\ref{6}), отримуємо
\begin{equation}\label{16}
X(t)=S(\Phi_{1}(t),\Phi_{2}(t))S^{-1}\exp \left \{ \frac{t}{T}S(\ln(Y(T)I_{0}))S^{-1} \right \} .
\end{equation}
Оскільки
\begin{equation}\label{17}
S(\ln(Y(T)I_{0}))S^{-1}=\ln(SY(T)S^{-1}SI_{0}S^{-1})=\ln(X(T)I),
\end{equation}
де
\begin{equation}\label{18}
I=SI_{0}S^{-1},
\end{equation}
то (\ref{16}) набуває вигляду необхідного нам представлення (ІІ) за
умови, що
\begin{equation}\label{19}
H=\frac{1}{T}\ln(X(T)I),
\end{equation}
\begin{equation}\label{20}
\Phi(t)=S(\Phi_{1}(t), \Phi_{2}(t))S^{-1}.
\end{equation}
Враховуючи (\ref{18}), отримуємо
\begin{equation}
I^{2}=E,\qquad I_{1}=S \left( \begin{array}{cc} E_1 & 0 \\ 0 & 0
\end{array} \right)S^{-1}, \qquad I_{2}=S \left( \begin{array}{cc} 0 & 0 \\ 0 &
E_2 \end{array} \right)S^{-1}, \nonumber
\end{equation}
де $E_{1}$ та $E_{2}$ --- одиничні матриці відповідних розмірів.

На підставі виразів для $I_{1}$ та $I_{2}$ одержуємо рівності
\begin{equation}\label{21}
\Phi(t)I_{1}=S(\Phi_{1}(t),0)S^{-1},\qquad\Phi(t)I_{2}=S(0,\Phi_{2}(t))S^{-1}
\end{equation}
для всіх $t\in\mathbb{R}$. Із (\ref{21}) з урахуванням властивостей
матриць $\Phi_{1}(t)$ та $\Phi_{2}(t),$ випливають співвідношення
\begin{equation}
\Phi(t+T)I_{1}=\Phi(t)I_{1},\qquad
\Phi(t+T)I_{2}=-\Phi(t)I_{2}\nonumber
\end{equation}
для всіх $t\in\mathbb{R}$, що й завершує доведення твердження 2
теореми  3 у розглядуваному випадку .

Нехай $d=n.$ Тоді замість (\ref{5}) отримуємо рівність
\begin{equation}
X(T)=SBS^{-1},\nonumber
\end{equation}
замість (\ref{8}) --- рівність
\begin{equation}
Y(t+kT)=(-1)^{k}Y(t)B_1^k,\nonumber
\end{equation}
замість (\ref{13}) --- рівність
\begin{equation}
Y(t)=\Phi_{2}(t) \exp \left \{ \frac{t}{T}\ln B_{1} \right \}
\nonumber
\end{equation}
та умову
\begin{equation}
\Phi_{2}(t+T)=-\Phi_{2}(t)\nonumber
\end{equation}
для всіх $t\in\mathbb{R}$.

Покладемо
\begin{equation}
I_{0}=-E\nonumber
\end{equation}
та  з двох останніх формул отримаємо рівність (ІІ) вигляду
\begin{equation}
Y(t)=\Phi_{2}(t) \exp \left \{ \frac{t}{T}\ln (-Y(t)) \right \}, \nonumber
\end{equation}
 де
\begin{equation}
H=\frac{1}{T}\ln (-X(T)), \qquad H \in \textbf{M}_{n}\mathbb{R},
\nonumber
\end{equation}
\begin{equation}
\Phi(t)=S\Phi_{2}(t)S^{-1}, \qquad  \Phi(t+T)=-\Phi(t),\qquad
\Phi(t) \in \textbf{M}_{n}\mathbb{R} \nonumber
\end{equation}
 для всіх $t\in\mathbb{R}$. Це й завершує доведення твердження 2
 теореми 3.

Перейдемо до доведення тверджень 3 теореми  3 та виокремимо з них
два граничні випадки: перший, коли матриця $X(t)$ має дійсний
логарифм, і другий, коли всі власні значення матриці $X(t)$ від'ємні
та серед їх елементарних дільників немає однакових.

У першому випадку твердження 3 теореми  3 випливають із формул
Г.~Флоке (ІІ), (ІІІ), згідно з якими заміна змінних
\begin{equation}
x=\Phi(t)z\nonumber
\end{equation}
зводить диференціальне рівняння (І) до диференціального рівняння
\begin{equation}
\frac{dz}{dt}=Hz\nonumber
\end{equation}
та забезпечує наявність вказаних в теоремі 3 властивостей матриць $H$ та
$\Phi(t)$.

У другому випадку твердження 3 теореми  3 є тривіальними:
заміна змінних
\begin{equation}
x=z\nonumber
\end{equation}
зводить диференціальне рівняння (І) до диференціального рівняння з тією ж матрицею коефіцієнтів
\begin{equation}
\frac{dz}{dt}=P(t)z.\nonumber
\end{equation}

Якщо пов'язати виокремлені вище граничні випадки з представленням
матриці $X(T)$ через її дійсну канонічну форму (\ref{5}), то першому
з них відповідає випадок, коли
\begin{equation}
X(T)=SAS^{-1},\nonumber
\end{equation}
другому --- випадок, коли
\begin{equation}
X(T)=SBS^{-1}.\nonumber
\end{equation}
Отже,  неграничним випадком тверджень 3 теореми 3 залишається
 лише той, який визначається умовами
\begin{equation}
A\in \textbf{M}_{n-d}(\mathbb{R}),\qquad B\in
\textbf{M}_{d}(\mathbb{R}), \qquad n>d>1.\nonumber
\end{equation}
Нехай ці умови виконуються. Тоді з доведення тверджень 2 теореми
 3 випливає, що матриця $Y(t)$, пов'язана з матрицею $X(t)$
формулою (\ref{6}), має вигляд (\ref{13}). Якщо позначити
\begin{equation}
U(t)=\Phi_{1}(t),\qquad V(t)=\Phi_{2}(t),\qquad H_{1}=\frac{\ln
A}{T}\,,\qquad H_{2}=\frac{\ln B_{1}}{T}, \nonumber
\end{equation}
то (\ref{13})  набуде вигляду
\begin{equation}\label{22}
Y(t)=(U(t),V(t))\left( \begin{array}{cc} e^{H_1 t} & 0 \\ 0 & e^{H_2 t} \end{array} \right).
\end{equation}
Із (\ref{22}) випливає, що
\begin{equation}\label{23}
Y(t)\left( \begin{array}{cc} E_1  \\ 0  \end{array} \right)=(U(t),V(t))\left( \begin{array}{cc} e^{H_1 t}  \\ 0  \end{array} \right)=U(t)e^{H_{1}t}\,,
\end{equation}
де $E_{1}$ --- $(n-d)$-вимірна одинична матриця.

Продиференціювавши рівність (\ref{23}) з урахуванням першої з формул
(\ref{6}),  отримуємо, що
\begin{equation}
S^{-1}P(t)SY(t)\left( \begin{array}{cc} E_1 \\ 0 \end{array} \right)=S^{-1}P(t)SU(t)e^{H_{1}t}=\frac{dU(t)}{dt}e^{H_{1}t}+U(t)H_{1}e^{H_{1}t}.
\nonumber
\end{equation}
Отже,
\begin{equation}\label{24}
\frac{dU(t)}{dt}+U(t)H_{1}=S^{-1}P(t)SU(t)
\end{equation}
для всіх $t\in\mathbb{R}$.

Матриця $Y(t)$ є фундаментальною матрицею розв'язків диференціального рівняння
\begin{equation}\label{25}
\frac{dy}{dt}=S^{-1}P(t)Sy.
\end{equation}

Нехай $W(t)\in \textbf{M}_{n \:d}(\mathbb{R})$ для всіх
$t\in\mathbb{R}$ та є неперервно диференційовною на $\mathbb{R}$,
періодичною  з періодом $T$ і такою, що
\begin{equation}\label{26}
\det(U(t),W(t))\ne 0
\end{equation}
для всіх $t\in\mathbb{R}$.

Існування такої матриці випливає з теореми про квазіперіодичний
базис в $\mathbb{R}^{n}$ роботи [\,4\,].

Зробимо в диференціальному рівнянні (\ref{25}) заміну змінних за формулою
\begin{equation}\label{27}
y=U(t)y_{1}+W(t)y_{2}\,
\end{equation}
та на підставі рівності (\ref{24}) отримаємо диференціальне рівняння
\begin{equation}\label{28}
U(t)\left(\frac{dy_{1}}{dt}-H_{1}y_{1}\right)+W(t)\frac{dy_{2}}{dt}=\left(S^{-1}P(t)SW(t)
-\frac{dW(t)}{dt}\right)y_{2}\,.
\end{equation}
Розв'язуючи це рівняння з допомогою матриці
\begin{equation}\label{29}
\left( \begin{array}{cc} L_1(t) \\ L_2(t) \end{array} \right),
\end{equation}
оберненої до матриці $(U(t),W(t))$, отримуємо для $\displaystyle
\frac{dy_{1}}{dt}$ та $\displaystyle \frac{dy_{2}}{dt}$ систему
диференціальних рівнянь
\begin{equation}\label{30}
\frac{dy_{1}}{dt}=H_{1}y_{1}+L_{1}(t)\left(S^{-1}P(t)SW(t)-\frac{dW(t)}{dt}\right)y_{2}\,,
\end{equation}
\begin{equation}\label{31}
\frac{dy_{2}}{dt}=L_{2}(t)\left(S^{-1}P(t)SW(t)-\frac{dW(t)}{dt}\right)y_{2}.
\end{equation}

Оскільки матриця коефіцієнтів системи (\ref{30}), (\ref{31}) має блочно-трикутний вигляд,  фундаментальною матрицею розв'язків цієї системи є матриця
\begin{equation}\label{32}
\left( \begin{array}{cc} e^{H_1 t} & Y_1(t) \\ 0 & Y_2(t) \end{array} \right),
\end{equation}
другий стовпчик якої утворений розв'язками системи диференціальних рівнянь
(\ref{30}), (\ref{31}) із заданим початковим значенням
$y_{1}=Y_{1}(0),\,\,y_{2}=Y_{2}(0)$ таким, що
\begin{equation}\label{33}
\det\,Y_{2}(0)\ne 0.
\end{equation}
Врахувуючи заміну змінних (\ref{27}), фундаментальною матрицею розв'язків рівняння (\ref{25}) є матриця
\begin{equation}\label{34}
(U(t),W(t))\left( \begin{array}{cc} e^{H_1 t} & Y_1(t) \\ 0 & Y_2(t) \end{array} \right).
\end{equation}
Крім того, формула (\ref{22}) також визначає  фундаментальну матрицю
розв'язків  рівняння (\ref{25}). Згідно з теорією лінійних
диференціальних рівнянь між цими двома фундаментальними матрицями
розв'язків існує залежність така, що
\begin{equation}\label{35}
(U(t),V(t))\left( \begin{array}{cc} e^{H_1 t} & 0 \\ 0 & e^{H_2 t} \end{array} \right)C=(U(t),W(t))\left( \begin{array}{cc} e^{H_1 t} & Y_1(t) \\ 0 & Y_2(t) \end{array} \right)
\end{equation}
для всіх $t\in\mathbb{R}$, де $C$ --- невироджена стала матриця. Підставивши в (\ref{35}) значення $t=0$, отримуємо для визначення матриці $C$ лінійне алгебраїчне рівняння
\begin{equation}\label{36}
(U(0),V(0))C=(U(0),W(0))\left( \begin{array}{cc} E_1 & Y_1(0) \\ 0 &
Y_2(0) \end{array} \right).
\end{equation}
Домноживши (\ref{36}) на матрицю $\left( \begin{array}{cc} L_1(0) \\
L_2(0) \end{array} \right),\,$  одержимо рівність
\begin{equation}\label{37}
\left( \begin{array}{cc} E_1 & L_1(0)V(0) \\ 0 & L_2(0)V(0) \end{array} \right)
C = \left( \begin{array}{cc} E_1 & Y_1(0) \\ 0 & Y_2(0) \end{array} \right),
\end{equation}
з якої випливає, по-перше, що
\begin{equation}\label{38}
\det(L_{2}(0),V(0))\ne 0,
\end{equation}
по-друге, що при
\begin{equation}\label{39}
Y_{1}(0)=L_{1}(0)V(0),\qquad Y_{2}(0)=L_{2}(0)V(0)
\end{equation}
одержимо
\begin{equation}\label{40}
C=E.
\end{equation}
Таким чином, визначивши за початковим значенням (\ref{39}) розв'язки
$\left( \begin{array}{cc}  Y_1(t) \\  Y_2(t) \end{array} \right)$ системи диференціальних
рівнянь (\ref{30}), (\ref{31}), отримаємо із (\ref{35}), (\ref{40}) рівність
\begin{equation}\label{41}
(U(t),V(t))\left( \begin{array}{cc} e^{H_1 t} & 0 \\ 0 & e^{H_2 t} \end{array} \right)=(U(t),W(t))\left( \begin{array}{cc} e^{H_1 t} & Y_1(t) \\ 0 & Y_2(t) \end{array} \right)
\end{equation}
для всіх $t\in\mathbb{R}$.

Помноживши (\ref{41}) на матрицю (\ref{29}), прийдемо до рівності
\begin{equation} \left( \begin{array}{cc} E_1 & L_1(t)V(t) \\ 0 & L_2(t)V(t)\end{array} \right) \left( \begin{array}{cc} e^{H_1 t} & 0 \\ 0 & e^{H_2 t} \end{array} \right)=\left( \begin{array}{cc} e^{H_1 t} & Y_1(t) \\ 0 & Y_2(t) \end{array} \right). \nonumber \end{equation}
Отже,
\begin{equation}\label{42}
Y_{1}(t)=L_{1}(t)V(t)e^{H_{2}t},
\end{equation}
\begin{equation}\label{43}
Y_{2}(t)=L_{2}(t)V(t)e^{H_{2}t}
\end{equation}
для всіх $t\in\mathbb{R}$.
Оскільки матриця $Y_{2}(t)$ є невиродженою, то із (\ref{43}) знаходимо значення $e^{H_{2}t}$ і, підставивши його  у рівність (\ref{42}), одержимо, що
\begin{equation}\label{44}
Y_{1}(t)=L_{1}(t)V(t)(L_{2}(t)V(t))^{-1}Y_{2}(t)
\end{equation}
для всіх $t\in\mathbb{R}$.

Запишемо систему диференціальних рівнянь (\ref{30}), (\ref{31}) у вигляді системи
\begin{equation}\label{45}
\frac{dy_{1}}{dt}=H_{1}y_{1}+R_{1}(t)y_{2},
\end{equation}
\begin{equation}\label{46}
\frac{dy_{2}}{dt}=G(t)y_{2},
\end{equation}
ввівши позначення
\begin{equation}
R_{1}(t)=L_{1}(t)\left(S^{-1}P(t)SW(t)-\frac{dW(t)}{dt}\right),\nonumber
\end{equation}
\begin{equation}
G(t)=L_{2}(t)\left(S^{-1}P(t)SW(t)-\frac{dW(t)}{dt}\right).\nonumber
\end{equation}
Запишемо також рівність (\ref{44}), використовуючи матрицю
\begin{equation}\label{47}
F(t)=L_{1}(t)V(t) (L_{2}(t)V(t))^{-1},
\end{equation}
у вигляді
\begin{equation}\label{48}
Y_{1}(t)=F(t)Y_{2}(t).
\end{equation}
Продиференціювавши рівність (\ref{48}) та врахувавши те, що матриця
$\left( \begin{array}{cc} Y_1(t) \\  Y_2(t) \end{array} \right)$ є
блоком фундаментальної матриці (\ref{32}) розв'язків системи
диференціальних рівнянь (\ref{30}), (\ref{31}), а отже, і системи
(\ref{45}), (\ref{46}), отримуємо співвідношення
\begin{equation}\label{49}
\frac{dF(t)}{dt}+F(t)G(t)=H_{1}F(t)+R_{1}(t)
\end{equation}
для всіх $t\in\mathbb{R}$.
Нарешті, зробимо заміну змінних
\begin{equation}\label{50}
y_{1}=z_{1}+F(t)z_{2},\qquad y_{2}=z_{2}
\end{equation}
та одержимо замість системи диференціальних рівнянь (\ref{45}), (\ref{46}) систему
\begin{equation}
\frac{dz_{1}(t)}{dt}+\frac{dF(t)}{dt}z_{2}+F(t)G(t)z_{2}=H_{1}z_{1}+H_{1}F(t)z_{2}+
R_{1}(t)z_{2},\nonumber
\end{equation}
\begin{equation}
\frac{dz_{2}(t)}{dt}=G(t)z_{2},\nonumber
\end{equation}
яка внаслідок (\ref{49}) набуде вигляду
\begin{equation}\label{51}
\frac{dz_{1}(t)}{dt}=H_{1}z_{1},\qquad
\frac{dz_{2}(t)}{dt}=G(t)z_{2}.
\end{equation}
Оскільки друге з рівнянь системи (\ref{51}) збігається з точністю до позначень із рівнянням (\ref{46}), то матриця $Y_{2}(t)$ є фундаментальною матрицею розв'язків другого з рівнянь системи (\ref{51}). Тоді згідно з формулою (\ref{43}) матриця
\begin{equation}\label{52}
L_{2}(t)V(t)e^{H_{2}t}(L_{2}(0)V(0))^{-1}
\end{equation}
є фундаментальною матрицею розв'язків другого з рівнянь системи
(\ref{51}), яка при $t=0$ є  одиничною матрицею. Отже, матриця
\begin{equation}\label{53}
(L_{2}(T)V(T))e^{H_{2}T}(L_{2}(0)V(0))^{-1}
\end{equation}
є матрицею монодромії другого з рівнянь системи (\ref{51}).

Згідно з  визначеннями матриця $L_{2}(t)$ є періодичною з періодом
$T$, матриця $V(t)$ задовольняє умову
\begin{equation}\label{54}
V(t+T)=-V(t),
\end{equation}
а матриця $H_{2}$ має вигляд
\begin{equation}\label{54}
H_{2}=\frac{1}{T}\ln (-B).\nonumber
\end{equation}
Враховуючи наведені властивості матриць $L_{2}(t)$, $V(t)$ та
$H_{2}$, приходимо до висновку, що матриця (\ref{53}) має вигляд
\begin{equation}
(-L_{2}(0)V(0))(-B)(L_{2}(0)V(0))^{-1}=L_{2}(0)V(0)B(L_{2}(0)V(0))^{-1}.\nonumber
\end{equation}
Отже, згідно з викладеним вище та  визначенням матриці $B$
отримуємо, що множина власних значень матриці (\ref{53}) є або
 множиною всіх від'ємних власних значень матриці $X(T)$, або
 її підмножиною.

Розглянемо матрицю $F(t)$. Із визначення цієї матриці за формулою
(\ref{47}) і того, що матриці $L_{1}(t)$ та $L_{2}(t)$ періодичні  з
періодом $T$, а матриця  $V(t)$ задовольняє умову (\ref{54}),
випливає, що
\begin{equation}
F(t+T)=(-L_{1}(t)V(t))(-L_{2}(t)V(t))^{-1}=F(t)\nonumber
\end{equation}
для всіх $t\in\mathbb{R}$.

Отже, матриця $F(t)$ є періодичною  з періодом $T$.

Для завершення доведення тверджень 3 теореми  3 залишилось взяти до
уваги те, що заміна змінних
\begin{equation}\label{55}
x=Sy
\end{equation}
перетворює диференціальне рівняння (І) у диференціальне рівняння (\ref{25}), тому
суперпозиція замін (\ref{55}), (\ref{27}), (\ref{25}) перетворює диференціальне
 рівняння (І) у систему диференціальних рівнянь (\ref{51}) з визначеними теоремою  3
 властивостями як заміни змінних, так і самих диференціальних рівнянь системи (\ref{51}).

Зробимо деякі зауваження до тверджень 2 та 3 теореми  3.

Перше з них стосується формули (VI), яка визначає матрицю $H$. Із
доведення теореми 3 випливає, що $H$ не завжди визначається
однозначно. Ця неоднозначність пов'язана з умовою розбиття дійсної
канонічної форми матриці $X(T)$ на блоки $A$ та $B,$ згідно з якою
матриця $B$ може бути як блоком жорданової форми матриці $X(T),$
утвореним всіма її жордановими клітинами, що відповідають її
від'ємним власним значенням, так і блоком цієї форми, утвореним із
визначеного вище блоку вилученням довільного числа пар однакових
жорданових клітин.

Друге зауваження стосується мінімально можливої величини порядку
другого з диференціальних рівнянь системи (VIII). Із доведення
теореми 3 випливає, що цей порядок також пов'язаний з умовами
розбиття дійсної канонічної форми матриці $X(T)$ на блоки $A$ та $B$
і дорівнює мінімально можливому порядку матриці $B$ такого розбиття.
Отже, як випливає з першого зауваження, мінімально можлива величина
порядку другого з рівнянь системи (VIII) дорівнює порядку матриці,
яка утворена із жорданової форми матриці $X(T)$ вилученням з неї
всіх жорданових клітин, що відповідають невід'ємним власним
значенням матриці $X(T),$ та максимально можливого парного числа
однакових жорданових клітин цієї ж матриці, які відповідають її
від'ємним власним значенням.

Крім того, відзначимо, що згідно з доведенням теореми  3 матриця $B$
є жордановою формою матриці монодромії $Z_2(T)$  другого з
 рівнянь системи  (VIII) з усіма наслідками
 щодо властивостей фундаментальної матриці розв'язків
$Z_2(t)$ цього рівняння, які звідси випливають.

Нарешті, зауважимо, що на підставі тверджень 3 теореми 3
та теореми 2 диференціальне рівняння (I) має інваріантні многовиди
\begin{equation}
K^{n-d}(t)=\{ x\in \mathbb{R}^n:\, U(t)L_1(t)x=x \}, \nonumber
\end{equation}
\begin{equation}
K^{d}(t)=\{ x\in \mathbb{R}^n:\, V(t)L_2(t)x=x \}, \nonumber
\end{equation}  періодичні  з періодом $T$:
\begin{equation}
K^{\nu}(t+T)=K^{\nu}(t), \qquad \nu \in \{(n-d) \vee d\} \nonumber
\end{equation} для всіх $t\in \mathbb{R} .$ При цьому (I) на
$K^{n-d}(t)$ еквівалентне першому, а на $K^{d}(t)$
--- другому з диференціальних рівнянь системи (VIII).

\textbf{Наслідок.} \emph{Фундаментальна матриця розв'язків
диференціального рівняння} (І) $X(t)$ \emph {має представлення}
$$ X(t)=\Phi (t) e^{Ht}\Phi ^+(0),  \eqno (\textrm{IX}) $$
\emph{де}
$$ H=\frac{1}{T}\ln \left( \begin{array}{cc} X(T) & 0 \\
0 & Z(T) \end{array} \right),   \eqno (\textrm{X}) $$ \emph{$H \in
\mathbf{M}_m(\mathbb{R}), \; Z(T)$
--- матриця монодромії звуження} (І) \emph {на його періодичний інваріантний
многовид $\,K^d(t), \;\Phi(t)$ --- періодична з періодом $T$
матриця, що задовольняє рівняння }
$$ \frac{d \Phi}{dt}+\Phi H=P(t)\Phi , \eqno (\textrm{XI}) $$
$\Phi(t)\in \mathbf{M}_{n\:m}(\mathbb{R}) $ \emph{для всіх $t
\in\mathbb{R} ,\; \Phi^+(0)$ --- матриця, псевдообернена для
матриці} $\Phi(0),\; m=n+d, \;   n\geq d \geq 0$.

Дійсно, згідно з останнім із зауважень диференціальне рівняння (І)
має періодичний інваріантний многовид $K^d(t),$ на якому рівняння
(І) еквівалентне другому з диференціальних рівнянь системи (VIII).
Розглянемо систему диференціальних рівнянь
\begin{equation}
\label{56}
\frac{dx}{dt}=P(t)x,\qquad   \frac{dz}{dt}=G(t)z,
\end{equation}
утворену з рівняння (І) та другого з рівнянь системи (VIII). Для її
матриці монодромії
\begin{equation}
\label{57} \left( \begin{array}{cc} X(T) & 0 \\
0 & Z(T) \end{array} \right)
\end{equation}
згідно з доведенням тверджень 3 теореми  3 дійсною
канонічною формою є матриця
\begin{equation}
\label{58} \left( \begin{array}{ccc} A & 0 & 0\\
0 & B & 0 \\ 0 & 0 & B  \end{array} \right),
\end{equation}
де $A$ та $B$ --- блоки розбиття дійсної канонічної форми  матриці
$X(T)$  такі, що матриця $A$ має дійсний логарифм. Оскільки матриця
\begin{equation}
\label{59} \left( \begin{array}{cc} B & 0 \\
0 & B \end{array} \right)
\end{equation}
утворена попарно однаковими жордановими клітинами, то й вона має
дійсний логарифм. Отже, логарифм матриці (\ref{58}) можна вибрати
дійсним. Тоді дійсним можна вибрати логарифм матриці (\ref{57}) і
визначити матрицю $H$ згідно з фор\-мулою~ (X) так, щоб вона
задовольняла умову $H \in \textbf{M}_{n+d}(\mathbb{R}),$ де $d$~---
порядок матриці $B.$ Застосувавши до фундаментальної матриці
розв'язків системи диференціальних рівнянь (\ref{56})
\begin{equation}
\left( \begin{array}{cc} X(t) & 0 \\
0 & Z(t) \end{array} \right) \nonumber \end{equation} формулу
Г.~Флоке (II), отримуємо, що
\begin{equation}
\label{60}
\left( \begin{array}{cc} X(t) & 0 \\
0 & Z(t) \end{array} \right) =
\left( \begin{array}{cc} \Phi_1(t)  \\
\Phi_2(t)  \end{array} \right)e^{Ht},
\end{equation}
де $H$ --- матриця (X) із простору $\textbf{M}_{n+d}(\mathbb{R}), \:
\Phi_1(t)$ та $\Phi_2(t)$
 --- періодичні  з періодом $T$ матриці,
$\Phi_1(t) \in \textbf{M}_{n \: n+d}(\mathbb{R})$ та $\Phi_2(t) \in
\textbf{M}_{d \: n+d}(\mathbb{R})$ для всіх $t\in\mathbb{R}$.
Диференціюючи рівність (\ref{60}), одержимо для матриці
$\Phi(t)=\left( \begin{array}{cc} \Phi_1(t)  \\
\Phi_2(t)  \end{array} \right) $ матричне диференціальне рівняння
\begin{equation}
\frac{d\Phi}{dt}+ \Phi H=\left( \begin{array}{cc} P(t) & 0 \\
0 & G(t) \end{array} \right)\Phi, \nonumber \end{equation} з якого
випливає, що матриця $\Phi_1(t)$ задовольняє диференціальне рівняння
(XI). Нарешті, помножимо рівність (\ref{60}) на матрицю
$\Phi_1^+(0)$, псевдообернену до матриці $\Phi_1(0)$, та отримаємо
рівність
\begin{equation} X(t)=\Phi_1(t)e^{Ht}\Phi_1^+(0), \nonumber \end{equation}
що збігається з (IX) з точністю до позначень.

\begin{center}\textbf{5. Два застосування отриманих результатів.}\end{center}

\textbf{ {1}.} \; Нехай $x \in \mathbb{R}^n, \: P(t)$ --- неперервна
періодична  з періодом $T$ матриця, $P(t)\in
\textbf{M}_{n}(\mathbb{R})$ для всіх $t\in\mathbb{R},\: X(t, x)$ ---
функція змінних $t\in\mathbb{R}$ та $x \in \mathbb{R}^n ,$ яка
набуває значень в $\mathbb{R}^n$, неперервна для всіх
$t\in\mathbb{R}$ та $x \in \mathbb{R}^n.$

Розглянемо диференціальне рівняння
\setcounter{equation}{0}
\begin{equation}
\label{1} \frac{dx}{dt}=P(t)x+ X(t, x).
\end{equation}
Нехай
\begin{equation}
\label{2} X(t, x)\equiv 0 ,
\end{equation}
тоді диференціальне рівняння (\ref{1}) має фундаментальну матрицю
розв'язків $X(t),$ для якої існує представлення
\begin{equation}
\label{3} X(t)=\Phi (t) e^{Ht} \Phi ^+(0) ,
\end{equation}
причому властивості матриць $\Phi (t)$ та $H$ визначені у наслідку
попереднього параграфа.

Використаємо формулу (\ref{3}) для зведення диференціального
рівняння (\ref{1}) до більш простого вигляду. Для цього зробимо в
(\ref{1}) заміну змінних, ввівши замість $x \in \mathbb{R}^n$ змінну
$y \in \mathbb{R}^m$ за формулою
\begin{equation}
\label{4} x=\Phi (t)y.
\end{equation}
Враховуючи те, що матриця $\Phi (t)$ є розв'язком диференціального
рівняння (XI), із (\ref{1}) та (\ref{4}) отримаємо рівність
\begin{equation}
\label{5} \Phi (t) \left(\frac{dy}{dt}- Hy \right)=X(t, \Phi (t)y).
\end{equation}
Подамо  цю рівність  у вигляді
\begin{equation}
\label{6} \frac{dy}{dt}- Hy= \Phi ^+(t) X(t, \Phi (t)y),
\end{equation}
де $\Phi ^+(t)$ --- матриця, псевдообернена для $\Phi (t)$ та має з
$\Phi (t)$ однакову гладкість та період, зокрема, за $\Phi ^+(t)$
можна вибрати перший блок матриці $(\Phi_1 ^+(t),\: \Phi_2 ^+(t)),$
оберненої до матриці $\left( \begin{array}{cc} \Phi_1(t)  \\
\Phi_2(t)  \end{array} \right)$ з формули (\ref{60}) попереднього
параграфа.  Розв'яжемо (\ref{6}) відносно $\displaystyle
\frac{dy}{dt}- Hy$ та
 отримаємо,
\begin{equation}
\label{7} \frac{dy}{dt}= Hy + \Phi ^+(t) X(t, \Phi (t)y).
\end{equation}
Виділена лінійна частина рівняння (\ref{7}) має постійну матрицю
коефіцієнтів, а загальна  зберігає властивості відповідної частини
початкового рівняння (\ref{1}).

\vspace{6.0mm}
 \textbf{ {2}.} \; Розглянемо диференціальне
рівняння
\begin{equation}
\label{8} \frac{dx}{dt}=X(x) + X_1(t, x),
\end{equation}
де $X(x)$ --- неперервно диференційовна функція  змінної $x, \;
X_1(t, x)$
--- непе\-рер\-вна функція змінних $t$ та $x ,$ яка набуває значень в $\mathbb{R}^n$ для
всіх $t \in \mathbb{R}$ та $x \in \mathbb{R}^n, \;  n\geq 2.$

Нехай рівняння (\ref{8}) за умови
\begin{equation}
\label{9} X_1(t, x)\equiv 0
\end{equation}
має періодичний розв'язок з періодом $T$
\begin{equation}
\label{10} x(t)=\xi (\omega t),
\end{equation}
де $\xi(\varphi)$ --- періодична за змінною $\varphi$ з періодом
$2\pi$ функція, $\omega =\displaystyle \frac{2\pi}{T}$ --- частота
періодичного розв'язку.

Рівняння у варіаціях, що відповідає періодичному розв'язку
(\ref{10}) диференціального рівняння (\ref{8}) при виконанні умови
(\ref{9}), має вигляд
\begin{equation}
\label{11} \frac{d\delta \xi}{dt}= \frac{\partial X(\xi (\omega
t))}{\partial x}\,\delta \xi
\end{equation}
та розв'язок
\begin{equation}
\label{12} \delta \xi=\xi '(\omega t),
\end{equation}
де $'$ означає похідну за змінною $\varphi$.

Дійсно, згідно з визначенням
\begin{equation}
\label{13} \xi ' (\varphi)\omega=X(\xi (\varphi)),
\end{equation}
отже,
\begin{equation}
\label{14} \xi '' (\varphi)\omega=\frac{\partial X(\xi (\varphi
))}{\partial x}\, \xi '(\varphi)
\end{equation}
для всіх $\varphi \in \mathbb{R}.$ Підставивши в (\ref{14}) замість
$\varphi$  значення $\omega t,$ отримаємо тотожність
\begin{equation}
\label{15} \frac{d }{dt}\,\xi '(\omega t)=\frac{\partial X(\xi
(\omega t ))}{\partial x}\, \xi '(\omega t),
\end{equation}
яка й доводить необхідне твердження.

Нехай $B(\varphi)$ --- неперервно диференційовна, періодична
 з періодом $2\pi$ матриця, $B(\varphi) \in
\textbf{M}_{n\,n-1}(\mathbb{R})$ та задовольняє умову
\begin{equation}
\label{16} \det (\xi ' (\varphi), B(\varphi))\neq 0
\end{equation}
для всіх $\varphi \in \mathbb{R}.$

Заміна змінних
\begin{equation}
\label{17} \delta \xi= \xi '(\omega t)c+B(\omega t)g
\end{equation}
зводить рівняння у варіаціях (\ref{11}) до диференціального рівняння
\begin{equation}
\xi^{''}(\omega t)\omega c +  \xi^{'}(\omega t)\frac{dc}{dt}+B\:{'}(\omega t)\omega g+
B(\omega t)\frac{dg}{dt}
= \frac{\partial X(\xi(\omega t))}{\partial x}(\xi^{'}(\omega t)c +
B(\omega t)g) \nonumber \end{equation}
або, з урахуванням (\ref{14}), до рівняння
\begin{equation}
\xi^{'}(\omega t)\frac{dc}{dt}+B(\omega t)\frac{dg}{dt}
=\left(\frac{\partial X(\xi(\omega t))}{\partial x}B(\omega t)-
B\:{'}(\omega t)\omega\right)g. \nonumber
\end{equation}
Це рівняння, якщо його розв'язати відносно похідних
$\displaystyle\frac{dc}{dt}$ та $\displaystyle \frac{dg}{dt}$ з
допомогою матриці
\begin{equation}
\label{18} \left( \begin{array}{cc} (\xi ' (\omega t))^+ \\ B^+(\omega t)  \end{array} \right),
\end{equation}
оберненої до $(\xi^{'}(\omega t),B(\omega t))$,  зводить (\ref{11})
до системи диференціальних рівнянь
\begin{align} \label{19}
\frac{dc}{dt}=(\xi ' (\omega t))^+
\left(\frac{\partial X(\xi(\omega t))}{\partial x} B(\omega t)-
B\:{'}(\omega t)\omega\right)g, \nonumber \\
\\
\frac{dg}{dt}=B^+(\omega t)\left(\frac{\partial X(\xi(\omega
t))}{\partial x} B(\omega t)- B\:{'}(\omega t)\omega\right)g
\nonumber.
\end{align}

У відповідності до заміни змінних (\ref{17}) матриця монодромії
системи диференціальних рівнянь (\ref{19}) подібна до матриці
монодромії рівняння у варіаціях (\ref{11}), таким чином, власні
числа обох матриць монодромії збігаються.

Із системи (\ref{19}) випливає, що одне з власних чисел її матриці монодромії дорівнює 1, а всі інші є власними числами матриці монодромії другого з диференціальних рівнянь системи (\ref{19}). Отже, такими ж є і власні числа матриці монодромії рівняння у варіаціях (\ref{11}).

Позначимо матрицю коефіцієнтів другого з диференціальних рівнянь системи (\ref{19})
через $Q(\omega t)$, де $Q(\varphi)$ --- періодична з періодом $2\pi$ матриця, та розглянемо це диференціальне рівняння
\begin{equation}\label{20}
\frac{dg}{dt}=Q(\omega t)g.
\end{equation}
На підставі  наслідку попереднього параграфа фундаментальна матриця
розв'яз\-ків рівняння (\ref{20}) $G(t)$ має представлення
\begin{equation}\label{21}
G(t)=\Phi (\omega t)e^{Ht}\Phi ^{+}(0),
\end{equation}
де
\begin{equation}\label{22}
H=\frac{1}{T} \ln \left( \begin{array}{cc} G(T) & 0 \\ 0 & Z(T)
\end{array} \right) \in \textbf{M}_{m}(\mathbb{R}),
\end{equation}
$2(n-1)\ge m\ge (n-1)\,,\,\,Z(t)$ --- фундаментальна матриця
звуження диференціального рівняння (\ref{20}) на його періодичний
інваріантний многовид $K^{m-(n-1)}(t)$, $\Phi (\varphi)$ ---
періодична  з періодом $2\pi$ матриця, $\Phi (\varphi)\in
\textbf{M}_{n-1 \:m}(\mathbb{R})$ для всіх $\varphi \in\mathbb{R}$
та задовольняє диференціальне рівняння
\begin{equation}\label{23}
\frac{d\Phi}{d\varphi}\,\omega +\Phi H=Q(\varphi )\Phi ,
\end{equation}
$\Phi^{+}(0)$ --- матриця, псевдообернена для $\Phi(0)$.

Використаємо отримані вище результати для введення амплітудно-фазових координат в околі замкненої кривої
\begin{equation}\label{24}
x=\xi(\varphi ),\,\,\varphi \in\mathbb{R} \, ,
\end{equation}
та для зведення диференціального рівняння (\ref{8}) в околі цієї
кривої до амплі\-туд\-но-фазової системи диференціальних рівнянь,
яка є більш простою.

Для цього у рівнянні (\ref{8}) зробимо заміну змінних згідно з
формулою
\begin{equation}\label{25}
x=\xi(\varphi)+B(\varphi)g,
\end{equation}
де $B(\varphi)$ --- матриця, визначена раніше.

Тоді, враховуючи рівність (\ref{13}), отримаємо замість (\ref{8}) диференціальне рівняння
\begin{equation}
(\xi^{'}(\varphi )+B^{'}(\varphi )g)\left(\frac{d\varphi}{dt}-\omega\right)
+B(\varphi)\frac{dg}{dt}=\nonumber
\end{equation}
\begin{equation}\label{26}
=X(\xi(\varphi)+B(\varphi)g)-X(\xi(\varphi))-
B^{'}(\varphi)\omega g+
X_{1}(t,\xi(\varphi)+B(\varphi)g).
\end{equation}
Розв'яжемо рівняння (\ref{26}) відносно
$\displaystyle\frac{d\varphi}{dt}-\omega$ та
$\displaystyle\frac{dg}{dt}$ з допомогою матриці
\begin{equation}\label{27}
\left( \begin{array}{cc} L_1(\varphi, g)  \\ L_2(\varphi, g)  \end{array} \right),
\end{equation}
оберненої до $(\xi^{'}(\varphi )+B^{'}(\varphi )g,B(\varphi))$.

Вибравши достатньо мале значення   $\delta>0$, матрицю (\ref{27})
легко побудувати для всіх
\begin{equation}\label{28}
\varphi\in\mathbb{R},\qquad ||g||\le\delta
\end{equation}
за матрицею (\ref{18}),  поклавши
\begin{equation}\label{29}
\left( \begin{array}{cc} L_1(\varphi, 0)  \\ L_2(\varphi, 0)
\end{array} \right)= \left( \begin{array}{cc} (\xi^{'}(\varphi ))^+
\\ B^+(\varphi)  \end{array} \right) .
\end{equation}
 Тоді з (\ref{26}) отримаємо систему диференціальних рівнянь
\begin{equation}\label{30}
\frac{d\varphi}{dt}=\omega+L_{1}(\varphi,g)[X(\xi(\varphi)+B(\varphi)g)-X(\xi(\varphi))+
X_{1}(t,\xi(\varphi)+B(\varphi)g)],
\end{equation}
\begin{equation}\label{31}
\frac{dg}{dt}=L_{2}(\varphi,g)[X(\xi(\varphi)+B(\varphi)g)-X(\xi(\varphi))+
X_{1}(t,\xi(\varphi)+B(\varphi)g)].
\end{equation}
Запишемо диференціальне рівняння (\ref{31}) у вигляді
\begin{equation}\label{32}
\frac{dg}{dt}=B^{+}(\varphi)\frac{\partial X(\xi(\varphi))}{\partial x}B(\varphi)g
+G(\varphi,g)+L_{2}(\varphi,g)X_{1}(t,\xi(\varphi)+B(\varphi)g),
\end{equation}
де через $G(\varphi,g)$ позначено функцію
\begin{equation}
L_{2}(\varphi,g)(X(\xi(\varphi)+B(\varphi)g)-X(\xi(\varphi))-
\frac{\partial X(\xi(\varphi))}{\partial x}B(\varphi)g)+\nonumber
\end{equation}
\begin{equation}\label{33}
+(L_{2}(\varphi,g)-L_{2}(\varphi,0)) \frac{\partial
X(\xi(\varphi))}{\partial x}B(\varphi)g,
\end{equation}
що задовольняє умови
\begin{equation}\label{34}
G(\varphi,0)=0,\qquad \frac{\partial G(\varphi,0)}{\partial g}=0.
\end{equation}
Згідно з визначенням матриці $Q(\omega t)$ матриця коефіцієнтів
виділеної лінійної частини диференціального рівняння (\ref{32})
збігається з матрицею $Q(\varphi)$.  Отже, (\ref{32}) набуває
вигляду диференціального рівняння
\begin{equation}\label{35}
\frac{dg}{dt}=Q(\varphi)g+G(\varphi,g)+L_{2}(\varphi,g)X_{1}(t,\xi(\varphi)+B(\varphi)g).
\end{equation}
Нехай $\Phi(\varphi)$ та $H$ --- матриці, визначені з представлення
(\ref{21}) фундаментальної матриці розв'язків диференціального
рівняння (\ref{20}). За допомогою цих матриць зробимо перетворення
системи диференціальних рівнянь (\ref{31}), (\ref{32}), поклавши
\begin{equation}\label{36}
g=\Phi(\varphi)h.
\end{equation}
В результаті отримаємо замість (\ref{35}) диференціальне рівняння
\begin{equation}
\Phi^{'}(\varphi)h+\Phi(\varphi)\left(\frac{dh}{dt}-Hh\right)+\Phi(\varphi)Hh=
Q(\varphi)\Phi(\varphi)h+G(\varphi,\Phi(\varphi)h)+\nonumber
\end{equation}
\begin{equation}\label{37}
+L_{2}(\varphi,\Phi(\varphi)h)X_{1}(t,\xi(\varphi)+B(\varphi)\Phi(\varphi)h)
\end{equation}
або, враховуючи, що  $\Phi(\varphi)$ є розв'язком диференціального
рівняння (\ref{23}),  \\ ---~рівняння
\begin{equation}\label{38}
\Phi(\varphi)\left(\frac{dh}{dt}-Hh\right)=G(\varphi,\Phi(\varphi)h)
+L_{2}(\varphi,\Phi(\varphi)h)X_{1}(t,\xi(\varphi)+B(\varphi)\Phi(\varphi)h).
\end{equation}
Розв'яжемо (\ref{38}) відносно $\displaystyle\frac{dh}{dt}-Hh$ з
допомогою матриці $\Phi^{+}(\varphi)$, псевдооберненої до
$\Phi(\varphi)$, та отримаємо рівність
\begin{equation}\label{39}
\frac{dh}{dt}=Hh+\Phi^{+}(\varphi)[G(\varphi,\Phi(\varphi)h)+
L_{2}(\varphi,\Phi(\varphi)h)X_{1}(t,\xi(\varphi)+B(\varphi)\Phi(\varphi)h)].
\end{equation}
На підставі викладеного вище заміна змінних (\ref{36}) зводить
систему диференціальних рівнянь (\ref{31}), (\ref{32}) до системи
\begin{equation}\label{40}
\frac{d\varphi}{dt}=\omega +f(t,\varphi,\Phi(\varphi)h),
\end{equation}
\begin{equation}\label{41}
\frac{dh}{dt}=Hh+\Phi^{+}(\varphi)F(t,\varphi,\Phi(\varphi)h),
\end{equation}
де $H$ --- матриця вигляду (\ref{22}),
\begin{equation}
f(t,\varphi,g)=L_{1}(\varphi,g)[X(\xi(\varphi)+B(\varphi)g)-X(\xi(\varphi))+
X_{1}(t,\xi(\varphi)+B(\varphi)g)],\nonumber
\end{equation}
\begin{equation}
F(t,\varphi,g)=L_{2}(\varphi,g)\left[X(\xi(\varphi)+B(\varphi)g)-X(\xi(\varphi))
-\frac{\partial X(\xi(\varphi))}{\partial x}B(\varphi)g\right]+\nonumber
\end{equation}
\begin{equation}
+(L_{2}(\varphi,g)-L_{2}(\varphi,0))\frac{\partial X(\xi(\varphi))}{\partial x}B(\varphi)g.\nonumber
\end{equation}
Система диференціальних рівнянь (\ref{40}), (\ref{41}) і є шуканою.

Таким чином, початкове диференціальне рівняння (\ref{8}) в околі
замкненої кривої (\ref{24}) зводиться суперпозицією замін
(\ref{25}), (\ref{36}), а отже, заміною змінних
\begin{equation}
x=\xi(\varphi)+B(\varphi)\Phi(\varphi)h,\nonumber
\end{equation}
до системи диференціальних рівнянь (\ref{40}), (\ref{41}), де
функції $f(t,\,\varphi,\,g)$ та  \\ $F(t,\, \varphi,\,g)$ неперервні
за змінними $t,\, \varphi,\,g$ для
$t\in\mathbb{R},\,\varphi\in\mathbb{R},
\,g\in\mathbb{R}^{n-1},\,||g||\le\delta,$ набувають значень в
$\mathbb{R}$ та $\mathbb{R}^{n-1}$ відповідно, періодичні за змінною
$\varphi$ з періодом $2\pi$, матриця $\Phi(\varphi)\in
\mathbf{M}_{n-1 \:m}(\mathbb{R})$ для всіх $\varphi\in\mathbb{R}$ та
є періодичною  з періодом $2\pi$, матриця $H\in
\mathbf{M}_{m}(\mathbb{R}) ,$ її  власними значеннями є числа
\begin{equation}
\frac{1}{T}\ln \lambda_{j},\qquad j=\overline{1, \,n-1},\nonumber
\end{equation}
та їх $p$-кратні повторення, $1\ge p_{j}\ge 0,\,\,\sum\limits_{j=1}^{n-1}p_{j}=m-(n-1),\,\,
\\2(n-1)\ge m\ge (n-1),\;1$
 та $\lambda_{1},...,\lambda_{n-1}$ --- власні значення матриці монодромії рівняння у варіаціях
(\ref{11}).

Зведення диференціальних рівнянь, що розглядалися вище, до рівнянь
із сталою матрицею коефіцієнтів  їх виділеної лінійної частини є
суттєвим для подальшого дослідження цих рівнянь. Підтвердження даної
думки знаходимо, зокрема, в роботах [\,5,\,6\,], проте там проблема
такого зведення розв'язана лише частково.
%
%
%
%


















\vspace*{36.0mm}

\noindent 1. \textit{Самойленко А.М. } Про інваріантні многовиди
лінійних диференціальних рівнянь. -- Київ, 2009.
-- 10~с. -- (Препринт / НАН України. Ін-т математики; № 2009.7).\\
2. \textit{Floquet  G.} Sur les \'{e}quations diff\'{e}rentielles
lin\'{e}aires \`{a} coefficients p\'{e}riodiques
// Ann. \'{E}cole Norm. Sup. -- 1883. -- № 12.
 -- P.~47 -- 88.\\
3. \textit{Гантмахер Ф.Р.} Теория матриц. --- М.:  Наука, 1988. -- 552~c.\\
4. \textit{Самойленко А.М.} Элементы математической теории
многочастотных колебаний. Инвариантные торы. --- М.:  Наука, 1987. -- 304~c.\\
5.  \emph{Боголюбов Н. Н., Митропольский Ю. А.\/} \, Асимптотические
методы в теории нелинейных колебаний. --- М., Физматгиз, 1963.
--- Николай Николаевич Боголюбов: Собр. научн. тр. в 12 томах, т. III.
--- М.: Наука, 2005. --- 605~с. \\
6. \textit{Самойленко А.М., Рекке Л.}  Условия синхронизации одной
колебательной системы //  Укр. мат. журн. --  2005. -- \textbf{57},
№ 7. -- С.~922 -- 945.\\